\documentclass[letterpaper, paper,11pt]{AAS}		% for final proceedings (20-page limit)
 \pdfoutput=1
\usepackage{bm}
\usepackage{amsmath}
\usepackage{amssymb}
\usepackage{subcaption} %provides mai subfigure
\usepackage{overcite}
\usepackage{float}
\usepackage{footnpag}		      	% make footnote symbols restart on each page

\DeclareMathOperator{\diag}{diag}

\begin{document}

\title{A multirate variational approach to simulation and optimal control for flexible spacecraft}

\author{Yana Lishkova\thanks{DPhil Student, Department
of Engineering Science, University of Oxford, Parks Road, OX1 3PJ Oxford, UK},  
Sina Ober-Bl\"{o}baum\thanks{Professor, Department
of Mathematics, University of Paderborn, Warburger Str. 100, D-33098 Paderborn, Germany },
Mark Cannon\thanks{Associate Professor, Department
of Engineering Science, University of Oxford, Parks Road, OX1 3PJ Oxford, UK},
Sigrid Leyendecker \thanks{Prof., Chair of Applied Dynamics, University of Erlangen-Nuremberg, Immerwahrstrasse 1, D-91058 Erlangen, Germany}
}

\maketitle{}

\begin{abstract}
We propose an optimal control method for simultaneous slewing and vibration control of flexible spacecraft. Considering dynamics on different time scales, the optimal control problem is discretized on micro and macro time grids using a multirate variational approach. The description of the system and the necessary optimality conditions are derived through the discrete Lagrange-d'Alembert principle. The discrete problem retains the conservation properties of the continuous model and achieves high fidelity simulation at a reduced computational cost.  Simulation results for a single-axis rotational maneuver demonstrate vibration suppression and achieve the same accuracy as the single rate method at reduced computational cost.
\end{abstract}

\section{Introduction}
To provide a wide range of %/(communication, navigation and other) 
services across the globe satellites are often equipped with various flexible appendages such as solar panels, antennas  and mechanical manipulators. In recent years the demand for improvements in satellite functionality and efficiency has led to more lightweight and flexible structure designs operating under stringent performance, positioning and energy usage requirements. 
Unfortunately even the simplest attitude maneuvers are capable of exciting vibrations in the flexible structures causing loss of pointing accuracy, decreased performance %, incorrect positioning 
or even structural damage.\cite{Reading116,Reading433} Thus there is a need for a computationally efficient control method capable of maximizing the system's performance while complying with hard safety and operational %critical
constraints. 

A control methodology capable of guaranteeing efficient operation while respecting such constraints is optimal control. Its first application to the problem of flexible satellite control can be traced back as early as the 1970s. Initial developments included Markley's work on performance index selection considering a single flexible mode.\cite{Reading492} The work of Breakwell attempted the control of several modes using feedback control and verified the results experimentally.\cite{Reading433} Turner and Junkins also presented work including several modes of vibration and explored the difference in results obtained using a nonlinear model or its linearization.\cite{Reading494} This research was later extended by Turner and Chun for %employing 
a distributed control system.\cite{Reading158} In all these works and many subsequent ones, the satellite was modelled as a distributed parameter system and approximated using a spatial discretization technique. This leads to the formulation of coupled ordinary and partial differential equations of motion. Two of the most common spatial discretization techniques are the Assumed Modes Method (AMM) and the Finite Element Method (FEM). A detailed review of both and their implementation for flexible spacecraft models can be found in Junkins and Kim. \cite{Reading157} However, it must be noted that in both of these techniques a more accurate representation of the flexible dynamics requires the construction of higher-dimensional models, which can lead to a computationally intensive search for an optimal control law.

Based on this description of the system, an Optimal Control Problem (OCP) can be formulated %/The solution of the resulting optimal control problem can be achieved 
and solved using broadly two groups of methods often classified as direct and indirect.\cite{Reading498,Reading499} Indirect methods employ Pontryagin's Maximum Principle to derive the necessary optimality conditions which are then discretized and solved numerically. 
The method however requires the introduction and numerical integration of adjoint variables thus increasing the dimensionality of the system by a factor of two. Furthermore, indirect methods require indepth understanding of the problem and the method and good initial guesses of the states and/or the adjoints. Additionally, contraints and terminal conditions can be difficult to incorporate.\cite{Reading498} Direct methods on the other hand provide improved robustness to the choice of initial guess and readily allow for the modification of path-constraints and boundary conditions. These methods convert the optimal control problem into a finite-dimensional nonlinear optimization problem through finite-dimensional parametrization of the controls or of both the controls and the states. The problem can then be solved using standard Nonlinear Programming (NLP) solvers in which
the necessary optimality conditions %for %the discrete problem 
 are derived through the Karush–Kuhn–Tucker conditions. A "good" initial guess is still required to initialize the NLP solver as the solution can often converge to a non-global minimum. One of the main drawbacks of direct methods is the need to solve large NLP problems, which can require significant computational resources. \cite{Reading498,Reading502,Reading27}

To improve computational efficiency %address this issue
while maintaining high accuracy of the solution, we propose the use of a direct transcription method known as Multirate Discrete Mechanics and Optimal Control (Multirate DMOC), which was introduced by Gail et al. and Junge et
al.\cite{Reading8,Reading70} This method provides high fidelity solutions at reduced computational cost by reducing the number of optimization variables and equality constraints as well as providing a sparse structure for the constraint Jacobian.  \cite{Reading452,Reading8} In DMOC both the description of the mechanical system and the necessary optimality conditions are derived though a discrete version of the Lagrange-d'Alembert principle. The resulting structure-preserving time-stepping equations serve as equality constraints for the optimization problem and allow for a discrete OCP formulation, which inherits the conservation properties of the continuous-time model and provides a sparse structure for the Jacobian of the constraint function. Preservation of symplecticity and symmetries in the Lagrangian allow for accurate representation of energy and/or momenta of the system where classical numerical integration techniques such as the standard Euler or Runge-Kutta methods introduce numerical dissipation which can be problematic in fuel and energy optimization problems.\cite{Reading70, Reading71, Reading8,Reading27}

The main improvement in computational efficiency comes from the observation that rigid-flexible structures exhibit dynamics on different time scales associated with the motion of the rigid body and the appendages, respectively. Integrating the whole system with small time steps would ensure stable integration of the fast dynamics, but leads to a large computational effort. This can be adressed through the use of a multirate formulation of the forced variational integrator within DMOC. This Multirate DMOC formulation allows for the simulation of the slow and fast dynamics to be carried out on seperate macro and micro time domains, respectively.\cite{Reading8} 
The lower number of macro time nodes on which the slow generalized coordinates are computed reduces the number of unknowns in the optimization as well as the dimensionality of the constraints. This results in lower computational cost compared to the single rate method while achieving comparable %/similar} 
accuracy in the solution as both slow and fast dynamics are resolved with appropriate time steps.\cite{Reading8,Reading452} Further reduction in the computational cost and memory can be achieved by the exploitation of the resulting sparse structure %that DMOC provides 
in the Jacobians of the constraint and cost functions. 

The presence of dynamics on different time scales was also noted by Azadi et al.~where the system was separated using singular perturbation theory.\cite{Reading104} In this work, however, the rigid motion and the vibration are controlled using different methods. In comparison the multirate DMOC formulation allows for a unified control methodology, which could allow for better optimality and constraint handling capabilities, while achieving high fidelity simulation and control at a reduced computational cost. Additionally, this method does not require for the equations of motion of the two subsystems to be decoupled and for systems with several slow motions, the number of generalized coordinates discretized on the macro scale can straightforwardly be extended within DMOC to tailor the method to the specific application and obtain further reductions in computational time and resources. 

In this work we aim to demonstrate the advantages of Multirate DMOC as a simulation and optimal control method for flexible spacecraft. For this purpose we formulate an optimal control problem to perform a single-axis rotation while leaving the appendages quiescent at the end of the maneuver. A general linear model of the spacecraft is chosen to allow for the validation of the proposed method against the analytical solution and several tests are carried out to demonstrate the numerical convergence and conservation properties of both the multirate variational integrator and the full Multirate DMOC scheme. Improvements in the computational efficiency compared to the single rate standard formulation are thoroughly investigated and demonstrated.

\section{Mathematical model and problem formulation}

\begin{figure}[htb]
\centering
    \includegraphics[width=4in]{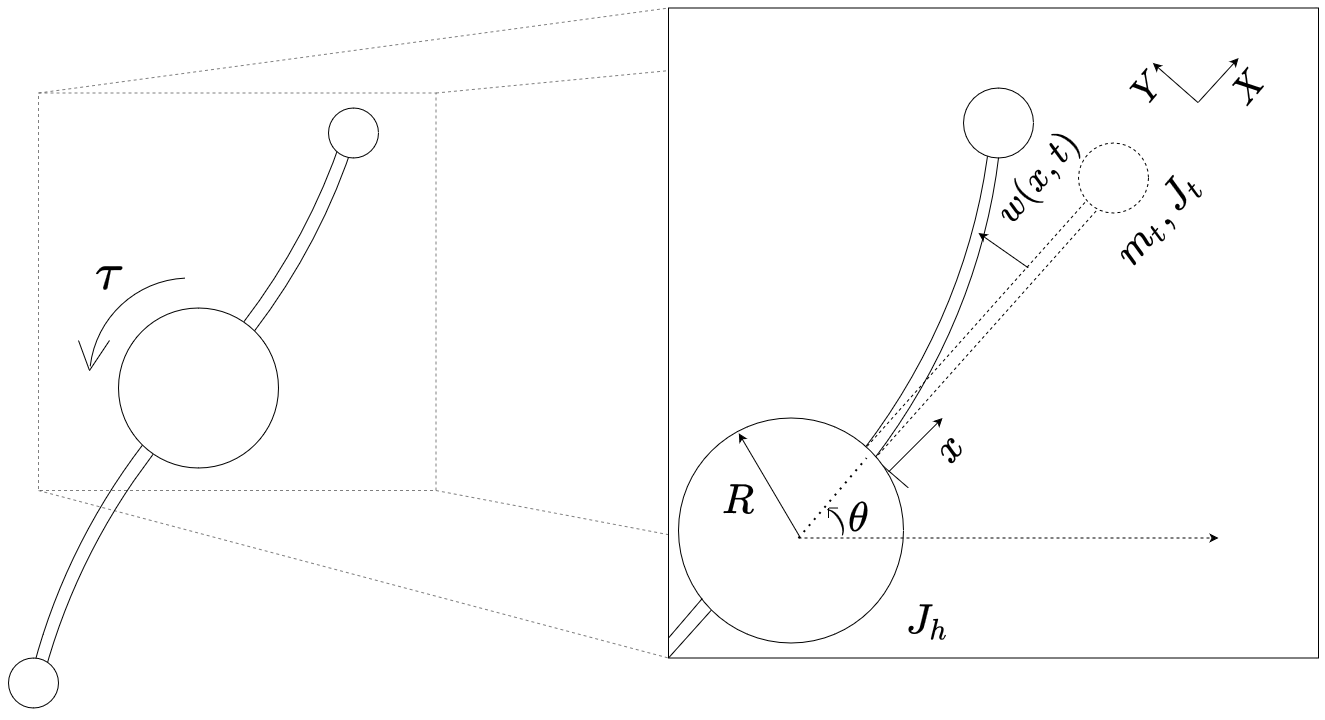}
    \caption{Spacecraft model during single-axis rotation with antisymmetric deformation}
    \label{fig:Sat1}
\end{figure}

To represent a wide range of flexible spacecraft we model the system as a rigid hub, representing the main body, with two rigidly-attached flexible appendages as shown in Figure~\ref{fig:Sat1}.  
A single-axis rotation maneuver is considered under the action of torque control input $\tau(t)$ at the hub. It is assumed to excite the flexible appendages antisymmetrically thus causing no shift in the center-of-mass of the spacecraft and resulting in the same deflection profiles in both appendages. % denoted as $w(x,t)$.  
These elastic transverse deformations are restricted to the plane X-Y that is perpendicular to the axis of rotation as shown in Figure~\ref{fig:Sat1}. The appendages are modelled as cantilever Euler-Bernoulli beams with attached tip-masses, thus the effects of shear deformation and sectional rotary inertia are neglected. 
The appendages are assumed to be identical in geometry and material properties and no structural damping is considered.  %(116).
\subsection{Spatial discretization}

To obtain a finite-dimensional model of the system we employ the Assumed Modes spatial discretization Method (AMM).\cite{Reading157} This method expresses the transverse elastic deflection of the beam $w(x,t)$ as a truncated series of products of assumed spatial mode shapes $\phi_{j}$ and time varying modal amplitudes  $\eta_{j}$ as follows
\begin{equation}
	\label{eq:w}
	w(x,t)=\sum_{j=1}^{N} \phi_{j}(x)\:\eta_{j}(t)
\end{equation} 
Here $N$ is the number of modes retained in the approximation, $x\in[0,L]$ denotes the position along the beam and $L$ is the length of the beam.
Together with the angle of hub rotation $\theta$, the amplitudes ($\eta_1(t),\ldots,\eta_{N}(t)$) become the set of the generalized coordinates describing the configuration of the system in time. \cite{Reading157} For the context of rigid-flexible satellites a widely used set of mode shapes is
\begin{equation}
	\label{eq:phi}
	\phi_j(x)= 1-\textrm{cos}\Big(\frac{j\pi x}{L}\Big)+\frac{1}{2}(-1)^{(j+1)}\Big(\frac{j\pi x}{L}\Big)^2
\end{equation} 
This set of functions satisfies the physical and geometrical boundary conditions for a clamped-free beam, but has been shown to produce satisfactory results for models which include tip-masses. \cite{Reading157} 
 
\subsection{The Lagrangian}
Detailed formulation of the linearized equations of motions can be found in Junkins and Kim and Turner and Chun. \cite{Reading157, Reading158}  
Based on these derivations the Lagrangian of the system can be expressed as follows:
\begin{equation}
	\label{eq:L0} \mathcal{L}(\underline{\xi},\:\underline{\dot{\xi}}) =\frac{1}{2}\:\underline{\dot{\xi}}^{T}\:\mathbf{M}\:\underline{\dot{\xi}}-\frac{1}{2}\:\underline{\xi}^{T}\:\mathbf{K}\:\underline{\xi}
\end{equation}
where
\begin{equation*}
	\label{eq:Lmatrices}
	\underline{\xi} = \begin{bmatrix}
	\:\theta\: \\
	\:\underline{\eta}\:
\end{bmatrix}, \; \;	   
	\mathbf{M}=  \begin{bmatrix}
         M_{\theta \theta}& M_{\theta \eta}^{T} \\
         M_{\theta \eta} & M_{\eta \eta}  
        \end{bmatrix},  \;\; 
       \mathbf{K}= \begin{bmatrix}
         0 & 0  \\
        0 &  K_{\eta \eta}
        \end{bmatrix}         
\end{equation*}

\noindent and the time-dependence has been dropped to simplify the notation. Here $\underline{\eta}\in R^{N\times1}$ is a vector of generalized coordinates defined as $\underline{\eta} = [ \eta_{1}, \: \eta_{2}, \: \ldots, \:  \eta_{N}]^{T}$ and $\underline{\xi}$ represents the configuration vector of the entire system. The elements of the mass and stiffness matrices are defined as
\begin{equation}
	\label{eq:a}
M_{\theta \theta} = J_h +2\:\Big[J_t+m_t(R+L)^2+\int_0^L \rho A (R+x)^2dx\Big]
\end{equation}
\begin{equation}
	\label{eq:ba}
[M_{\theta \eta}]_i = 2\:m_t\:(R+L)\phi_i(L)+2\:J_t\,\phi_i'(L)+2\int_0^L \rho A (R+x)\phi_i(x)dx
\end{equation}
\begin{equation}
	\label{eq:bv}
[M_{\eta \eta}]_{i,j} = 2\:m_t\:\phi_i(L)\phi_j(L)+2\:J_t\,\phi_i'(L)\phi_j'(L)+2\int_0^L \rho A \phi_i(x)\phi_j(x)dx
\end{equation}
\begin{equation}
	\label{eq:bd}
[K_{\eta \eta}]_{i,j} = 2\:\int_0^L EI\phi''_i(x)\phi''_j(x)dx
\end{equation}

 \noindent for $M_{\theta \theta}\in R^{1\times1} $, $ M_{\theta \eta}\in R^{N\times1} $, $ M_{ \eta \eta}\in R^{N\times N} $ , $ K_{\eta \eta} \in R^{N\times N}$. Here $[\cdot]'$ and $[.]''$ denote the first and second derivative with respect to $x$, $[\cdot]_i$ denotes the i-th element of the respective vector and $[\cdot]_{i,j}$ the $(i,j)$-th element of the respective matrix. $J_h$ and $R$ denote the rotary inertia  and radius of the hub. $m_t$ and $J_t$ represent the mass and the rotary inertia of each of the tip-masses and $\rho$, $A$, $EI$ are respectively the density, the cross-sectional area and the flexural rigidity of the beams.
%These expressions 

The virtual work in this example can be shown to be
\begin{equation}
	\label{eq:WirtWork}
\delta W=\mathfrak{\underline{f}}\cdot\delta\underline{\xi} = \tau \:\delta \theta
\end{equation}
where $\mathfrak{\underline{f}}$ is the vector of generalized forces. Using the Lagrange-d'Alembert Principle, which requires that
\begin{equation}
	\label{eq:LdA1}
\delta \int_{t_0}^{t_f} \mathcal{L}(\underline{\xi},\:\underline{\dot{\xi}})\:dt+ \int_{t_0}^{t_f} (\,\mathfrak{\underline{f}}\cdot\delta\underline{\xi}\,)\:dt=0
\end{equation}

\noindent for all variations $\delta \underline{\xi}$ with $\delta \underline{\xi}(t_0) = \delta \underline{\xi}(t_f)=0$, one can then obtain the following equations of motion for this example spacecraft 
\begin{equation}
	\label{eq:EqMot}
\mathbf{M}\,\underline{\ddot{\xi}}+\mathbf{K}\,\underline{\xi} = \mathbf{D}\,\tau
\end{equation}
for $\mathbf{D} = [1,0,...,0]^{T}$. \cite{Reading157,Reading71}

\subsection{Transformation to modal coordinates} %TRANSFORMATION TO MODAL COORDINATES
Eq.~(\ref{eq:EqMot}) represents a set of $N+1$ coupled differential equations and to decouple them we introduce the Modal Coordinate Transformation.\cite{Reading252,Reading157} %They can be de/uncoupled using the natural modes of the system by the use of the modal coordinate transformation. 
Solving for the eigenvalue problem for $\mathbf{M}$ and $\mathbf{K}$, we obtain a set of eigenvalues $\lambda_{j}$ and normalized eigenvectors $\underline{e}_j$ as follows 

\begin{equation}
	\label{eq:MT4}
\left|\mathbf{K}-\lambda_j\,\mathbf{M}\right|=0,\; \;\mathbf{K}\: \underline{e}_{j}=\lambda_j\,\mathbf{M}\: \underline{e}_{j} , \; \; \; \lambda_{j}\leq\lambda_{j+1}   \; \; \;\textrm{for}\; \; \; j=1,2, \: ... \: ,N+1
\end{equation} 
 such that 
\begin{equation}
	\label{eq:MT1}
\; \; \; \mathbf{E}^{T}\mathbf{M}\mathbf{E} = \mathbf{I}, \; \; \; \mathbf{E}^{T}\mathbf{K}\mathbf{E} = \mathbf{\Lambda}
\end{equation}
where
\begin{equation*}
	\label{eq:MT5}
\mathbf{E}=[\underline{e}_{1}, \, \underline{e}_{2}, \, \ldots ,\, \underline{e}_{N+1}],\; \; \; \;  \mathbf{\Lambda} = \left[
\begin{array}{ccc} \lambda_{1} & \cdots & 0 \\
   \vdots & \ddots & \vdots \\
   0 & \cdots & \lambda_{N+1} 
\end{array}
\right] 
\end{equation*} 
 and $\mathbf{I}$ is the identity matrix. Physically  each $\lambda_{j}$ corresponds to the square of the natural frequency $w_{j}$ (rad/s) of the system under consideration and the corresponding eigenvector describes the corresponding mode shape. Thus introducing 
%/defining} 
the linear transformation
\begin{equation}
	\label{eq:MTransf}
\underline{\xi} = \mathbf{E}\:\underline{q}, \; \; \; \underline{\dot{\xi}} = \mathbf{E\:}\underline{\dot{q}}
\end{equation} 
the Lagrangian and the generalized force vector can be rewritten as
\begin{equation}
	\label{eq:LMT}
\mathcal{L} = \frac{1}{2}(\,\underline{\dot{q}}^{T}\underline{\dot{q}}-\underline{q}^{T}\mathbf{\Lambda}\,\underline{q}\,), \; \; \;\; \; \underline{f} = \mathbf{E}^{T}\mathbf{D}\:\tau
\end{equation}
%for $\mathbf{D} = [1,0,...,0]^{T}$. 
Applying the Lagrange-d'Alembert Principle as before, the equations of motion can be rewritten as 
\begin{equation}
	\label{eq:EqMotMT}
\underline{\ddot{q}}+\mathbf{\Lambda}\,\underline{q} = \mathbf{E}^{T}\mathbf{D}\,\tau
\end{equation}
This is a system of $N+1$ decoupled differential equations each describing a motion with natural frequency $w_{j}$.

\subsection{Problem formulation}

Having defined the model we proceed to formulate the optimal control problem with the objective of achieving a rest-to-rest single-axis rotation with $N$ quiescent modes at the end of the maneuver. %\textit{The performance index is chosen} 
The cost function is formulated as a sum of weighted quadratic functions of the modal amplitudes, the modal amplitude rates and the applied control and the overall OCP can be expressed as follows
\begin{subequations}
\begin{gather}
% \begin{equation}
	\label{eq:state0}
		J(\,x,u)= \int_{t_0}^{t_f}C(\,\underline{x}(t),u(t)\,)\:dt = \frac{1}{2}\int_{t_{0}}^{t_{f}}[\,\underline{x}(t)^T\mathbf{W}\underline{x}(t)+ u(t)^2\,]dt \\
% \end{equation}
 \notag \textrm{subject to \;\;\;\;\;\;\;\;\;\;\;\;\;\;\;\;\;\;\;\;\;\;\;\;\;\;\;\;\;\;\;\;\;\;\;\;\;\;\;\;\;\;\;\;\;\;\;\;\;\;\;\;\;\;\;\;\;\;\;\;\;\;\;\;\;\;\;\;\;\;\;\;\;\;\;\;\;\;\;\;\;\;\;\;\;\;\;\;\;\;\;\;\;\;\;\;\;\;\;\;\;\;\;\;\;\;\;\;\;\;\;\;\;\;\;\;\;\;\;\;\;\;\;\;\;\; } \\ 
%\begin{equation}
%	\label{eq:state2}
		\underline{\dot{x}}(t)=\mathbf{A}\underline{x}(t)+\mathbf{B}\tau(t) \\
%\end{equation} 
%\begin{equation}
%	\label{eq:state4}
		\underline{q}(t_{0}) =\mathbf{E}^{-1}\:\underline{\xi}_{\,t_0},\; \; \; \underline{\xi}_{\,t_0}= [0, ..., 0]^T  \\
%\end{equation} 
%\begin{equation}
%	\label{eq:state5}
		\; \;\; \;\underline{q}(t_{f}) = \mathbf{E}^{-1}\:\underline{\xi}_{\,t_f},\; \; \; \underline{\xi}_{\,t_f} =[\theta_{t_f},0,..., 0]^T \\
%\end{equation} 
%\begin{equation}
%	\label{eq:state7}
		\underline{\dot{q}}(t_{0}) = \underline{\dot{q}}(t_{f})=[0,..., 0]^T  \\
%\end{equation} 
	\notag \textrm{where \;\;\;\;\;\;\;\;\;\;\;\;\;\;\;\;\;\;\;\;\;\;\;\;\;\;\;\;\;\;\;\;\;\;\;\;\;\;\;\;\;\;\;\;\;\;\;\;\;\;\;\;\;\;\;\;\;\;\;\;\;\;\;\;\;\;\;\;\;\;\;\;\;\;\;\;\;\;\;\;\;\;\;\;\;\;\;\;\;\;\;\;\;\;\;\;\;\;\;\;\;\;\;\;\;\;\;\;\;\;\;\;\;\;\;\;\;\;\;\;\;\;\;\;\;\;\;\;\;\;\; }  \\
%\begin{equation*}
%	\label{eq:state8}
		\underline{x}(t) = \begin{bmatrix}
	\underline{q}(t) \\
	\underline{\dot{q}}(t)
\end{bmatrix}, \; \; u(t) = \tau(t), \; \;\mathbf{A} = \begin{bmatrix}
	0 & \mathbf{I}\, \\
	-\mathbf{\Lambda} & 0\,
\end{bmatrix},\; \;  \mathbf{B} = \begin{bmatrix}
	0 \\
	\mathbf{E}^{T}\mathbf{D}
\end{bmatrix} 
%\end{equation*} 
\end{gather}
\label{eq:ProbForm}
\end{subequations}
and $\mathbf{W}$ is taken to be the identity matrix for this example. This formulation of the cost function allows for the minimization of control effort, while penalizing trajectories with large modal deformations, which could lead to potential loss of pointing accuracy, degradation of performance or even structural damage. \cite{Reading116,Reading433} 

 \section{Simulation and optimal control using Multirate DMOC}
 \subsection{Multirate configuration description and discretization}
 %Based on the decoupled system description introduced previously, the system dynamics can further be \textcolor{blue}{seperated} ...
  
As shown previously the system can be decoupled into $N+1$ equations of motion. Depending on the natural frequency of each motion the system can further be separated into slow and fast subsystems as follows%can be described in terms of slow and fast subsystems %seperate the description of the system into a slow and fast subsystem as follows:
\begin{equation}
	\label{eq:two}
\underline{q}= \begin{bmatrix}
         \underline{q}^{s}  \\
        \underline{q}^{f} 
        \end{bmatrix}, \; \; \; \; \; \; \; \; 
        \underline{f}=
        \begin{bmatrix}
         \underline{f}^{s}  \\
        \underline{f}^{f}
        \end{bmatrix}=\mathbf{E}^T\mathbf{D}\: \tau 
        =\begin{bmatrix} \underline{Z}^s  \\
        \underline{Z}^f
        \end{bmatrix} \tau
\end{equation}

 \begin{equation}
	\label{eq:four}
 \mathcal{L} = \frac{1}{2}\Big(\:(\underline{\dot{q}}^{s})^{T}\underline{\dot{q}}^{s}-(\underline{q}^{s})^{T}\mathbf{\Lambda_{s}}\:\underline{q}^{s}\:\Big)+\frac{1}{2}\Big(\:(\underline{\dot{q}^{f}})^{T}\underline{\dot{q}}^{f}-(\underline{q}^{f})^{T}\mathbf{\Lambda_{f}}\:\underline{q}^{f}\:\Big)
\end{equation}
\vspace{2pt}

\noindent where $\mathbf{\Lambda_s} = \diag(\lambda_{1},\; \dots \;, \lambda_{r})$, $\mathbf{\Lambda_f} = \diag(\lambda_{r+1},\; \dots \;, \lambda_{N+1})$, $\; \underline{q}^s,\underline{\dot{q}}^s,\underline{f}^s,\underline{Z}^s \in R^{\:r\times1}$, $\; \underline{q}^f, \underline{\dot{q}}^f,\underline{f}^f, \underline{Z}^f \in R^{\:(N+1-r)\times 1}$ and the time dependence has been dropped to simplify the notation. The size of the slow subsystem $r$ is a free variable and %\textcolor{blue}{
expresses the number of modal coordinates, which will be treated as slow dynamics and thus discretized on a coarser time grid in the multirate approach. 

\begin{figure}[htb]
	\centering\includegraphics[width=\textwidth]{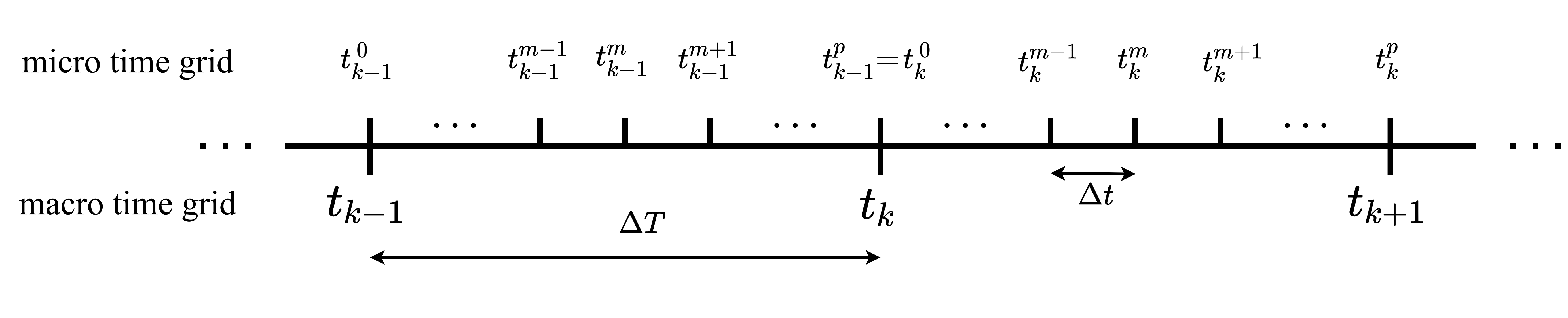}
	\caption{Micro and macro time grid schematic}
	\label{fig:Grid}
\end{figure}

In the discrete setting, we introduce two time grids: a macro time grid with  the macro time step $\Delta T$ and a micro time grid, which is obtained by subdividing each macro step into $p$ equally spaced micro time steps of size $\Delta t$ as depicted in Figure \ref{fig:Grid}. On the macro time grid we discretize the slow subsystem and on the micro time grid we discretize the fast dynamics. For this purpose define the discrete paths 

\begin{table}[H]
%\centering %% uncomment this if the "array" should be centered
\setlength\arraycolsep{10pt} % default value: 5pt
\def\arraystretch{2}
$\begin{array}{lll}           % "ll" sets up two left-aligned columns
 \;\; \;  \underline{q}_{\,d}^{s}(t_k) = \underline{q}_{\,k}^{s} & \textrm{for} \; \; t_k = t_0+ k\Delta T & \textrm{and} \:\: k=0,\ldots,n_s \\
 \;  \; \;\underline{q}_{\,d}^{f}(t_k^m) = \underline{q}_{\,k}^{f,m} & \textrm{for} \; \; t_k^m=t_0+k\Delta T+m\Delta t & \textrm{and} \;\;k=0,\ldots,n_s-1,\; m =0,\ldots,p \\
\end{array}$
\end{table} 

\noindent by approximating the  slow and the fast configuration variables using piecewise-linear polynomials as follows
 \begin{equation}
	\label{eq:five}
\underline{q}_{\,d}^{s}(t) = \underline{q}_{\,k}^{s}+\frac{\underline{q}_{\,k+1}^{s}-\underline{q}_{\,k}^{s}}{\Delta T}\:(t-t_k), \; \; \; \;  \; \; \; \;  \; \; \; \; \; \; \;  \;  \; \; \underline{\dot{q}}_{\,d}^{s}(t) = \frac{\underline{q}_{\,k+1}^{s}-\underline{q}_{\,k}^{s}}{\Delta T},  \; \; \; \; \; \; \;  \textrm{for} \; \;  t\in[t_k,t_{k+1}]
\end{equation}
 \begin{equation}
	\label{eq:seven}
\underline{q}_{\,d}^{f,m}(t) = \underline{q}_{\,k}^{f,m}+\frac{\underline{q}_{\,k}^{f,m+1}-\underline{q}_{\,k}^{f,m}}{\Delta t}\:(t-t_k^m), \; \; \; \; \underline{\dot{q}}_{\,d}^{f,m}(t) = \frac{\underline{q}_{\,k}^{f,m+1}-\underline{q}_{\,k}^{f,m}}{\Delta t}, \; \; \textrm{for} \; \; t\in[t_k^m,t_{k}^{m+1}]
\end{equation}
\vspace{2pt}

% \begin{equation}
\noindent where $n_s = t_f/\Delta T$ and denotes the number of macro time steps. On the other hand the control path $u = \tau$ is approximated by piecewise-constant values in $u_d=\tau_d$ which are defined at the midpoints of the micro time grid as follows
\begin{table}[!htbp]
%\centering %% uncomment this if the "array" should be centered
\setlength\arraycolsep{10pt} % default value: 5pt
\def\arraystretch{2}
$\begin{array}{ll}           % "ll" sets up two left-aligned columns
\; \; \;\;\; \; \;\; \; \;\;\; \; \;\;\; \; \;\; \; \;\;\; \;  \; \;\tau_d =\{ \{\tau_k^{m+1/2}\}^{p-1}_{m=0}\}^{n_s-1}_{k=0} & \textrm{where} \; \; \; \; \tau_k^{m+1/2} \approx \tau(t_k^{m+1/2}) \\
\end{array}$
\end{table}

\noindent Additionally we define $\underline{q}_{\,k}^f=\{\underline{q}_{\,k}^{f,m}\}^{p}_{m=0}$ and  $\tau_k=\{\tau_k^{m+1/2}\}^{p-1}_{m=0}$.

\vspace*{4pt}

\subsection{Discrete equations of motion}

The discrete equations of motion are derived following Gail et al.~and Ober-Bl{\"o}baum et al. \cite{Reading8,Reading71}. In standard direct transcription methods this is achieved through the direct discretization of the equations of motion, which we previously obtained using the Lagrange-d'Alembert principle (Eq.~(\ref{eq:EqMotMT})). In Multirate DMOC, however, the discretization is done one step earlier. Namely, the variation principle is discretized instead and the discrete equations of motion are derived from a discrete version of the above-mentioned Lagrange-d'Alembert principle {(Eq.~(\ref{eq:LdA1}))}.

First the discrete Lagrangian is defined as an approximation of the action integral for one macro time step

 \begin{equation}
	\label{eq:Lint}
\mathcal{L}_{k} = \mathcal{L}_d(\underline{q}_{\,k}^{s},\underline{q}_{\,k+1}^{s},\underline{q}_{\,k}^f) \approx \int_{t_k}^{t_{k+1}} \mathcal{L}(\underline{q}^s,\underline{q}^f,\underline{\dot{q}}^s,\underline{\dot{q}}^f) dt
\end{equation}
and thus the discrete action sum is defined as 
 \begin{equation}
	\label{eq:ActInt}
	S_d(\underline{q}_{\,d}^s,\underline{q}_{\,d}^{f}) = \sum^{n_s-1}_{k=0} \mathcal{L}_d(\underline{q}_{\,k}^s,\underline{q}_{\,k+1}^s,\underline{q}_{\,k}^f)
\end{equation}
Similarly on each macro step the virtual work of the external forces  is approximated as
%Similarly on each macro step we approximate the virtual work of the external forces  as:
 \begin{equation}
	\label{eq:ForceInt}
	\underline{f}_{\,k}^{s-}\cdot\delta\underline{q}_{\,k}^{s}+ \underline{f}_{\,k}^{s+}\cdot\delta\underline{q}_{\,k+1}^{s}+\sum^{p-1}_{m=0}(\underline{f}_{\,k}^{f,m-}\cdot\delta\underline{q}_{\,k}^{f,m}+ \underline{f}_{\,k}^{m+}\cdot\delta\underline{q}_{\,k}^{f,m+1}) \approx \int^{t_{k+1}}_{t_k}(\underline{f}^s\cdot\delta\underline{q}^s+\underline{f}^f\cdot\delta\underline{q}^f)dt
\end{equation}

\noindent where the control forces $\underline{f}^f$ and $\underline{f}^s$ have been discretized by discrete left and right forces as follows \cite{Reading57} 

\begin{equation}
 \underline{f}^{s\pm}_k=\underline{Z}^s \Big(\frac{\Delta t}{2}\sum_{m=0}^{p-1} \tau_k^{m+1/2} \Big) \; \; \; \; \textrm{for} \; \; \; \;  k=0,\ldots,n_s-1\; \; \; \;\; \; \; \;\; \; \; \;\; \; \; \;\; \; \; \;\; \; \; \;\; \; \; \;\;
\end{equation}
\begin{equation}
 \underline{f}^{f,m\pm}_k=\underline{Z}^f\Big(\frac{\Delta t}{2}\: \tau_k^{m+1/2} \Big)\; \; \; \; \; \; \; \textrm{for} \; \; \; \;  k=0,\ldots,n_s-1, \: m=0,\ldots,p-1
\end{equation}

Thus the discrete multirate Lagrange-d'Alembert principle can be expressed as
 \begin{equation}
	\label{eq:eight}
	\delta\sum^{n_s-1}_{k=0}(\mathcal{L}_k)+\sum^{n_s-1}_{k=0}\big[\, \underline{f}_{\,k}^{s-}\cdot\delta\underline{q}_{\,k}^s+\underline{f}_{\,k}^{s+}\cdot\delta\underline{q}_{\,k+1}^s+\sum^{p-1}_{m=0}(\:\underline{f}_{\,k}^{f,m-}\cdot\delta\underline{q}^{f,m}_{\,k}+\underline{f}_{\,k}^{f,m+}\cdot\delta\underline{q}_{\,k}^{f,m+1})\big]=0
\end{equation}

\noindent for all variations $\delta \underline{q}^s$ and $\delta \underline{q}^f$ vanishing at the end points ($\delta\underline{q}^s_{\,0} = \delta\underline{q}^s_{\,n_s}=\delta\underline{q}^{f,\,0}_{\,0}=\delta\underline{q}^{f,\,p}_{\,n_s-1}=0$). The condition for stationary curves $\underline{q}^s_{\,d},\;\underline{q}^f_{\,d}$ can be expressed as
\begin{subequations}
\begin{gather}
% \begin{equation}	
\label{eq:DiscrEqMot11}
 D_{\underline{q}_{\,k}^s}(\mathcal{L}_{k-1}+\mathcal{L}_{k})+\underline{f}_{\,k}^{s-}+\underline{f}_{\,k-1}^{s+}=0 \\
%\end{equation}
% \begin{equation}
	\label{eq:DiscrEqMot12}
D_{\underline{q}_{\,k}^{f,0}}(\mathcal{L}_{k-1}+\mathcal{L}_{k})+\underline{f}_{\,k}^{f,0-}+\underline{f}_{\,k-1}^{f,p-1+}=0 \\
%\end{equation}
% \begin{equation}
	\label{eq:DiscrEqMot14}
D_{\underline{q}_{\,k}^{f,m}}(\mathcal{L}_{k})+\underline{f}_{\,k}^{f,m-}+\underline{f}_{\,k}^{f,m-1+}=0 
%\end{equation}
\end{gather}
\label{eq:DiscrEqMot1}
\end{subequations}
%see 30 also
\noindent for $k=1,\,\ldots,n_s-1$ in Eq.~(\ref{eq:DiscrEqMot11}) and Eq.~(\ref{eq:DiscrEqMot12}) and $k=0,\,\ldots,n_s-1$ and $m=1,\,\ldots,p-1$ in Eq.~(\ref{eq:DiscrEqMot14}). Here $D_{[\cdot]}$ denotes the derivative with respect to the specified argument. These equations are known as the discrete forced multirate Euler-Lagrange equations and describe dynamics of the system under consideration in the discrete multirate setting. \cite{Reading8,Reading71}

\subsection{Discrete cost function and boundary conditions}

Similarly to the approximation of the Lagrangian and the virtual work, the cost functional (Eq. (\ref{eq:state0})) is approximated on each macro time step as
\begin{equation}
\label{eq:CostInt}
C_d(\underline{q}_{\,k}^s,\underline{q}_{\,k+1}^s,\underline{q}_{\,k}^f,\tau_k) = \int_{t_{k}}^{t_{k+1}}C(\underline{q}^s,\underline{q}^f,\underline{\dot{q}}^s,\underline{\dot{q}}^f,\tau)\:dt
\end{equation}

Boundary configuration and velocity conditions can be incorporated using the definitions of the fast and slow conjugate momentum, which are obtained using the Legendre transform as shown by Gail et al. and Leyendecker et al. \cite{Reading8,Reading9} 

\begin{table}[!htbp]
%\centering %% uncomment this if the "array" should be centered
\setlength\arraycolsep{10pt} % default value: 5pt
\def\arraystretch{2}
\begin{equation}
\label{eq:ConjMomentum}
\begin{array}{ll}           % "ll" sets up two left-aligned columns
\underline{p}_{\,k}^{s-} = -D_{\underline{q}^s_{\,k}}\:\mathcal{L}_d(\underline{q}_{\,k}^{s},\:\underline{q}_{\,k+1}^{s},\:\underline{q}_{\,k}^{f})-\underline{f}_{\,k}^{s-}  & \textrm{for}\;\; k=0,\ldots,n_s-1 \\
\underline{p}_{\,k}^{s+} = D_{\underline{q}^s_{\,k}}\:\mathcal{L}_d(\underline{q}_{\,k-1}^{s},\:\underline{q}_{\,k}^{s},\:\underline{q}_{\,k-1}^{f})+\underline{f}_{\,k-1}^{s+}        &\textrm{for}\;\; k=1,\ldots,n_s\\
\underline{p}_{\,k}^{f,0-} = -D_{\underline{q}^{f,0}_{\,k}}\:\mathcal{L}_d(\underline{q}_{\,k}^{s},\:\underline{q}_{\,k+1}^{s},\:\underline{q}_{\,k}^{f})-\underline{f}_{\,k}^{f,0-} & \textrm{for}\;\; k=0,\ldots,n_s-1\\
  \underline{p}_{\,k-1}^{f,p+} = D_{\underline{q}^{f,p}_{\,k-1}}\:\mathcal{L}_d(\underline{q}_{\,k-1}^{s},\:\underline{q}_{\,k}^{s},\:\underline{q}_{\,k-1}^{f})+\underline{f}_{\,k-1}^{f,p-1+}           &\textrm{for}\;\; k=1,\ldots,n_s \\
\end{array}
\end{equation}
\end{table}
\subsection{Discrete problem formulation}
In conclusion, the discrete constrained optimal control problem as defined by Multirate DMOC seeks to minimize the discrete cost functional
\begin{equation}
J_d(\underline{q}_{\,d}^s,\underline{q}_{\,d}^f,\tau_d)=\sum^{n_s-1}_{k=0}C_d(\underline{q}_{\,k}^s,\underline{q}_{\,k+1}^s,\underline{q}_{\,k}^f,\tau_k)
\end{equation}
subject to the discrete path constraints derived by Eq.~(\ref{eq:DiscrEqMot1}) and the boundary conditions enforced through Eq. (\ref{eq:ConjMomentum}). 

\section{Computational results}
Based on this formulation, multirate DMOC was implemented using MATLAB to simulate and optimally control the system. For this specific implementation we used the Midpoint Quadrature Rule to approximate the relevant integrals in Eq.~(\ref{eq:Lint}), Eq.~(\ref{eq:ForceInt}) and Eq.~(\ref{eq:CostInt}). \cite{Reading71} An example system is used for all simulations with the key geometrical and material parameters presented in Table~\ref{tab:Param} as suggested by Junkins and Kim. \cite{Reading157} %(158 does it in a text form-looks better). 
The model is constructed with five terms in the AMM approximation ($N\,$=$\,5$) and three degrees of freedom were treated as slow dynamics in the multirate formulation ($r\,$=$\,3$). Based on these specifications the natural frequencies of the system were obtained as $w_{1} = 0$ rad/s, $w_{2} = 6.454$ rad/s, $w_{3} = 52.41$, rad/s, $w_{4} = 1.607\times 10^{2}$ rad/s, $w_{5} = 3.381\times 10^{2}$ rad/s and $w_{N+1}=5.78\times 10^{2}$ rad/s resulting in dynamics on significantly different time-scales and allowing for thorough investigation of the advantages of the proposed multirate scheme.

\begin{table}[htbp]
	\fontsize{10}{10}\selectfont
    \caption{Structural parameters used for the simulations}
   \label{tab:Param}
        \centering 
   \begin{tabular}{l | l | c | l } % Column formatting, 
      \hline 
      Hub radius  & $R$ & 1.0 & ft\\
      Hub rotary inertia   & $J_h$ & 8.0 & slug-ft$^2$\\
      Tip mass   & $m_t$ & 0.156941 & slug\\
      Tip mass rotary inertia  & $J_t$ & 0.0018 & slug-ft$^2$\\              
      Beam length  & $L$ & 4.0 & ft\\   
      Beam height& $h$ & 6.0 & in.\\  
      Beam thickness& $t$ & 0.125 & in.\\  
      Beam linear density  & $\rho A$ & 0.0271875 & slug/ft\\
      Beam elastic modulus & $E$ & 0.1584$\times 10^{10}$ & lb$/$ft$^2$\\
      \hline
   \end{tabular}
\end{table}

\subsection{System simulation}

Before implementing the optimal control method, the proposed numerical integration scheme was validated against a reference analytic solution and its advantages were thoroughly investigated. This was achieved through a series of forward simulations in the absence of control for a time period of $t\in[0,t_{f}]$ starting from a rest position with initial deflection in the appendages of $\underline{\eta}(t_0)\,$=$[0.05,0.001,0.001,0.0001,0.0001]^T$.
%=$\,0.05\times10^{-(i-1)},i\,$=$\,1,\ldots,N$+$1$.}%

First, the numerical convergence of the method was investigated by computing the error of the simulated configuration with respect to the analytical continuous solution. The comparison was conducted on the macro time nodes for both the slow and fast variables ($\underline{q}^{s},\underline{q}^{f}$)  using the following error definition %\textcolor{red}{(30)}
\begin{equation}	
e_{\underline{x}}(\underline{x}_{\,d},\underline{x})= \max_{k\:=\:0,...,n_s} (\| \underline{x}_{\,k}-\underline{x}(t_k)\|_{\infty})
\label{eq:error}
\end{equation}

\begin{figure}[htb]
  \begin{subfigure}[b]{0.48\textwidth}
    \includegraphics[width=\textwidth]{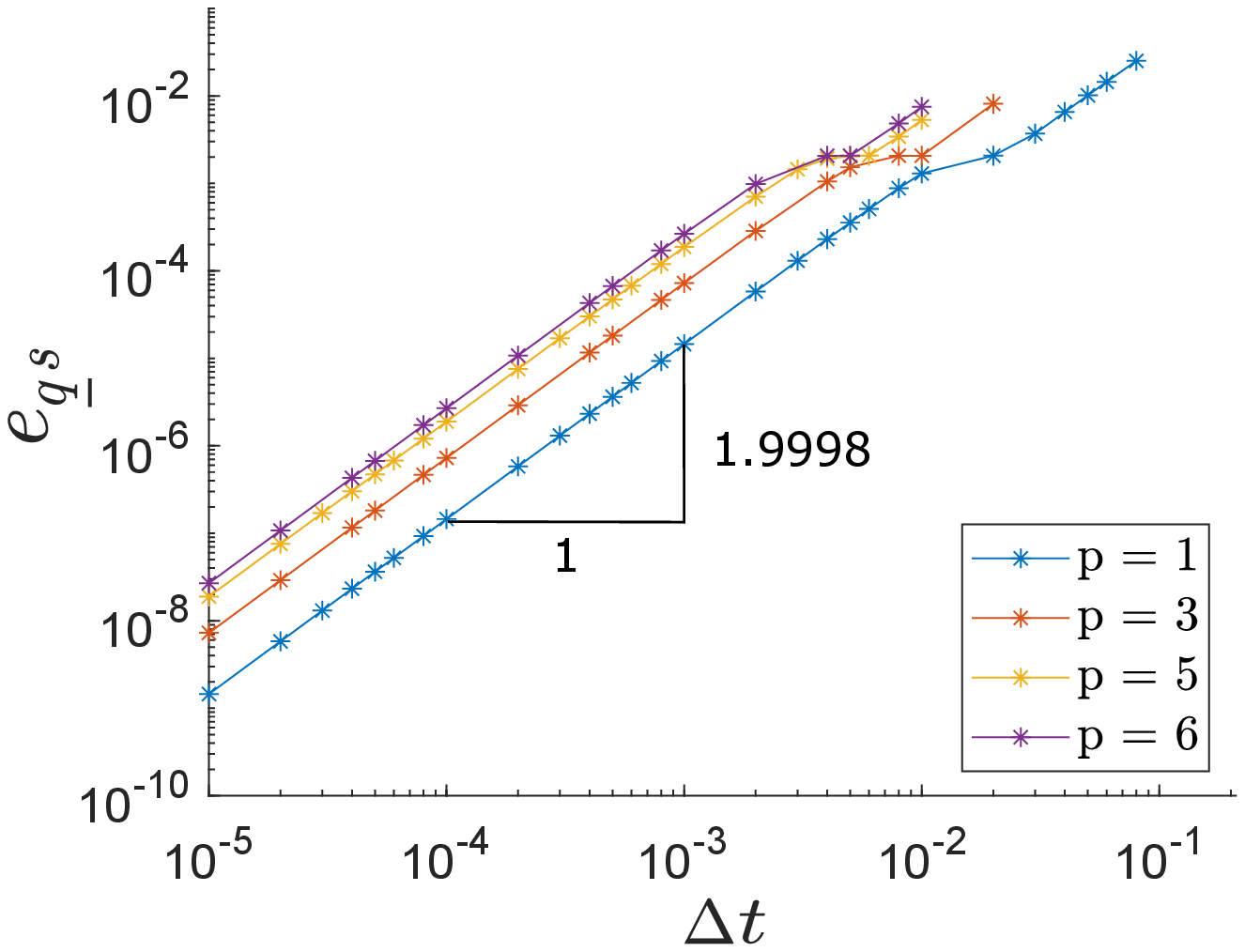}
    \caption{Error estimation for $\underline{q}^{s}$}
    \label{fig:Order1}
  \end{subfigure}
  \hfill
  \begin{subfigure}[b]{0.48\textwidth}
    \includegraphics[width=\textwidth]{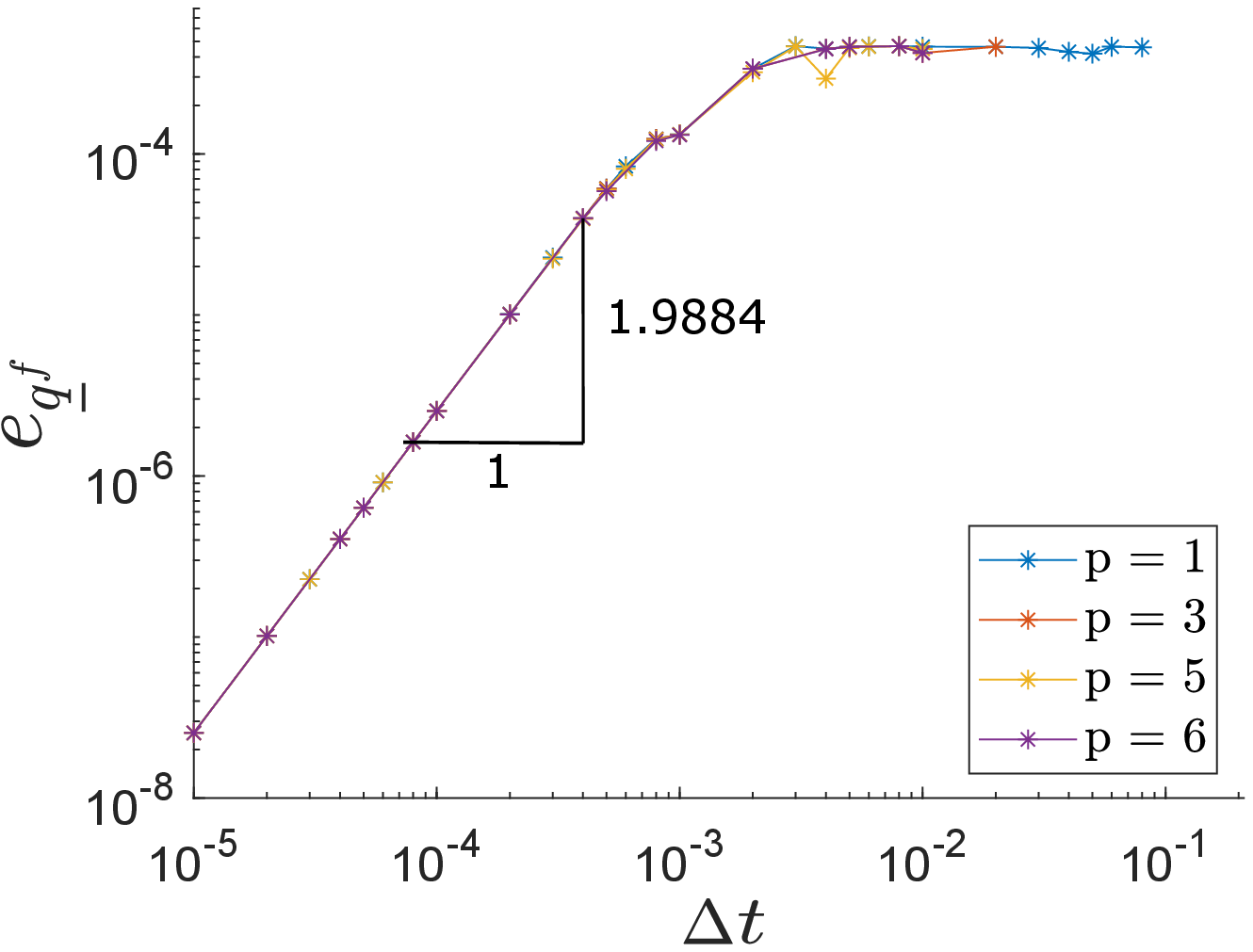}
    \caption{Error estimation for $\underline{q}^{f}$}
    \label{fig:Order2}
  \end{subfigure}
\caption{Numerical convergence study for the numerical integrator using macro-node error computation for $t_f=1.2s$}
\end{figure}

Figures \ref{fig:Order1} and \ref{fig:Order2} present the results for $t_f=1.2s$ and macro-micro proportionalities of $p\,$=$\,1,3$, $5,6$, where $p\,$=$\,1$ corresponds to a single rate solution. These figures demonstrate error convergence of approximately order 2, %(of the error) 
%for %all \textcolor{blue}{variables}  
as expected from the order of the chosen quadrature, thus validating the correct behaviour of the numerical integration scheme. We also note that an order reduction is observed in both figures for a small region of time-steps indicated by the presence of plateau, as expected for stiff systems (see Simeon).\cite{simeon1998order}
\begin{figure}[htb]
\begin{minipage}{0.48\textwidth}
	\centering
    \includegraphics[width=\textwidth]{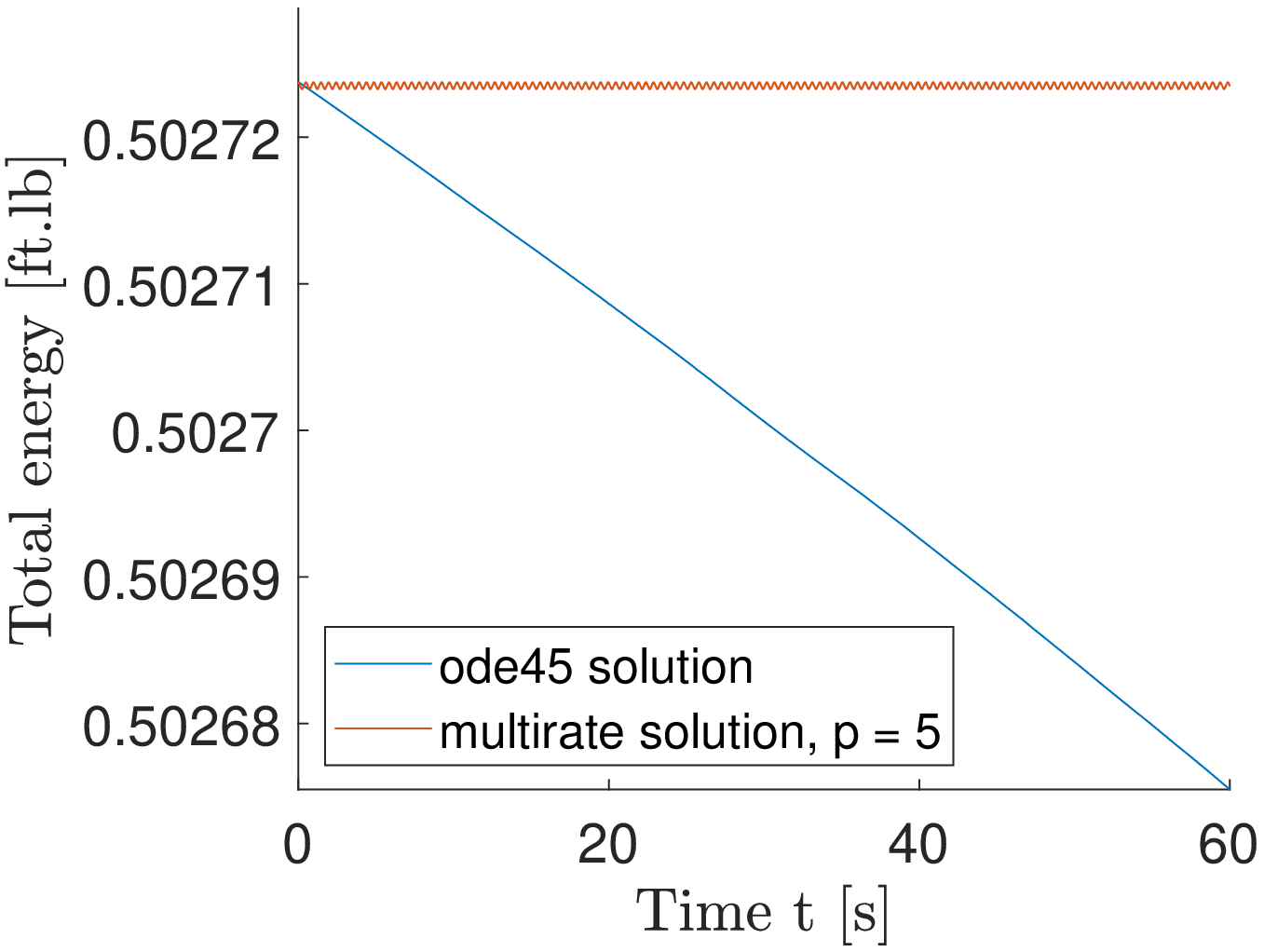}%energy60s
    \caption{Numerical dissipation of total energy for simulations with $\Delta t=10^{-4}$ and $t_f=60$s}
    \label{fig:TotEn}
\end{minipage}
  \hfill
\begin{minipage}{0.48\textwidth}
	\centering
    \includegraphics[width=\textwidth]{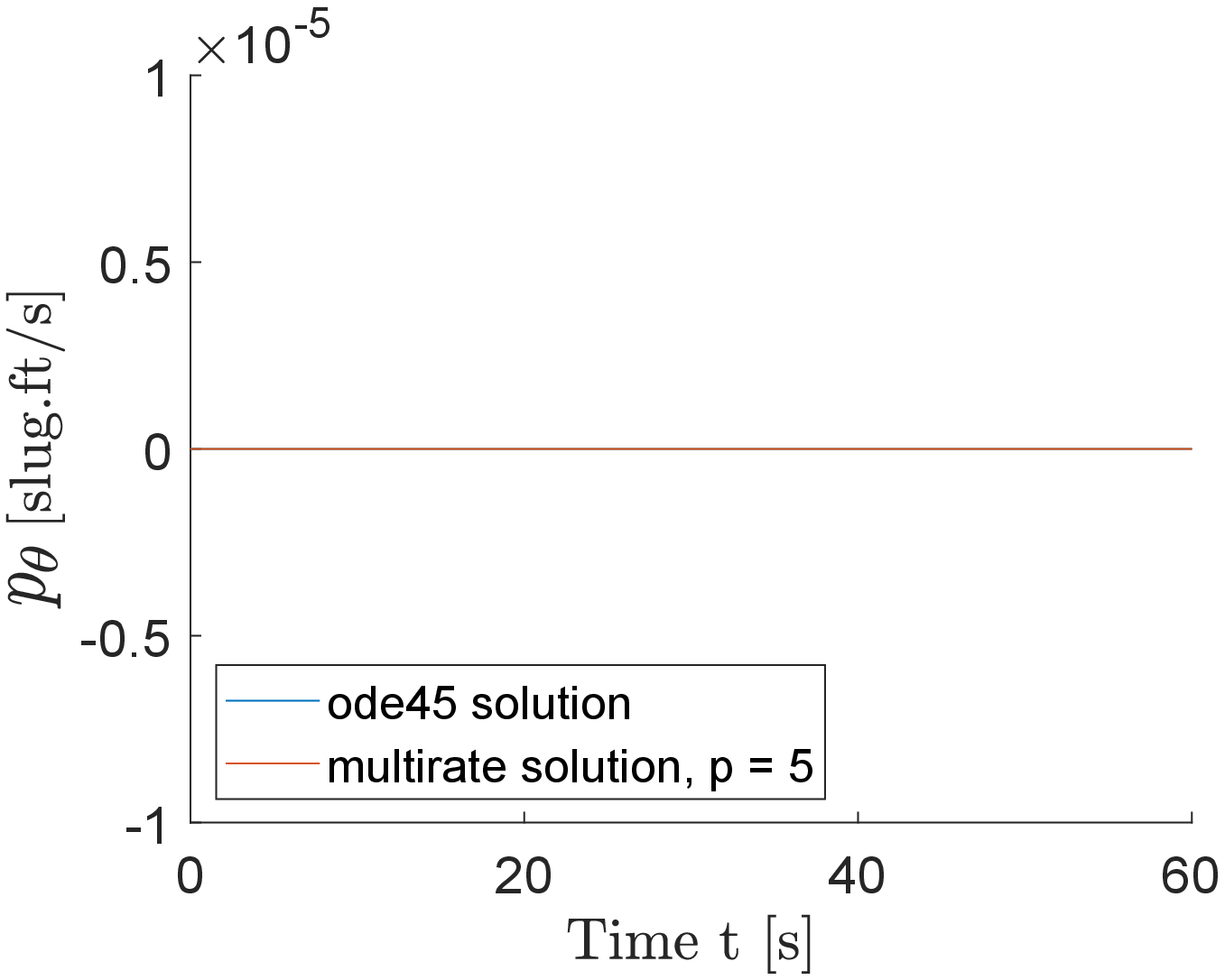}%momentum60s
    \caption{Momentum preservation for simulations with $\Delta t=10^{-4}$ and $t_f=60$s  }
    \label{fig:Moment}
\end{minipage}\hfill
\end{figure}

Figures~\ref{fig:TotEn} and \ref{fig:Moment} present the evolution of the total energy and the generalized momentum for $\theta$ defined as: 
\begin{equation}
p_{\theta} = \frac{\partial \mathcal{L}}{\partial \dot{\theta}}= M_{\theta \theta} \,\dot{\theta}+(M_{\theta \eta})^T \,\underline{\dot{\eta}}
\end{equation}
The results are obtained from simulations with step size of $\Delta t = 10^{-4}$ using both the variational scheme with $p\,$=$\,5$ and the \textit{ode45} MATLAB  numerical solver with \textit{RelTol} of $10^{-10}$. As the Lagrangian (Eq.~(\ref{eq:L0})) is not explicitly dependent on $t$ and $\theta$ and  no damping and external forces are considered in this test case, both the total energy and $p_{\theta}$ should remain constant in time according to Noether's theorem. \cite{Reading57} Figure~\ref{fig:Moment} demonstrates that both the \textit{ode45} implementation and the multirate solution successfully preserve $p_{\theta}$. In Figure~\ref{fig:TotEn}, however, it can be observed that the \textit{ode45} implementation introduces numerical dissipation, while the variational scheme successfully preserves the the total energy of the system up to small bounded fluctuations %oscillations
 known to be present for variational integrators (see West). \cite{Reading514} %, while \textit{ode45} implementation introduces to numerical dissipation.
Together, the two tests demonstrate the structure-preserving properties of the variational scheme.

\begin{figure}[!htb]
	\centering\includegraphics[width=0.48\textwidth]{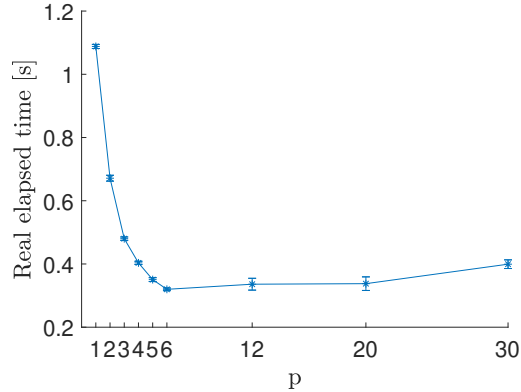} %
    \caption{Mean computational time with standard deviation versus $p$ for a constant
micro time step of $10^{-4}$ and $t_f = 1.2s$}
	\label{fig:CPUtimeIntegr2}
\end{figure}

\vspace*{-5pt}
Returning to Figures \ref{fig:Order1} and \ref{fig:Order2}, we can further observe that making the macro grid coarser by increasing $p$, increases the error in the computation of the slow variables slightly, but maintains the same accuracy in the approximation of the fast ones. This small sacrifice in the accuracy, however, enables great computational savings. To demonstrate the reduction in computational cost for %\textcolor{blue}{To demonstrate the computational efficiency 
the multirate approach %we ran a series of 
we performed a series of simulations for a number of macro-to-micro time step proportionalities $p$ while keeping the micro step constant at $10^{-4}$. Figure \ref{fig:CPUtimeIntegr2} presents the mean of the real elapsed time for the respective simulation obtained from 10 measurements using the \textit{tic-toc} MATLAB routine. A clear reduction in the required computational time can be seen by increasing $p$ up to an optimal value ($p\,$=$\,6$ in this instance). The existence of an optimum $p$-value is due to the specific implementation of the numerical integrator, which uses the Newton-Raphson Method to resolve the system one %\parfillskip=0pt\par}
%\noindent 
macro time step at a time. Increasing the number of micro time steps within this period increases the number of variables being obtained at each iteration of the method and thus results in slower computations. Nevertheless the results demonstrate that significant computational savings are possible for a wide range of $p$-values.

\subsection{Optimal control using Multirate DMOC}
\begin{figure}[htb]
\begin{minipage}{0.315\textwidth}
	\centering
    \includegraphics[width=\textwidth]{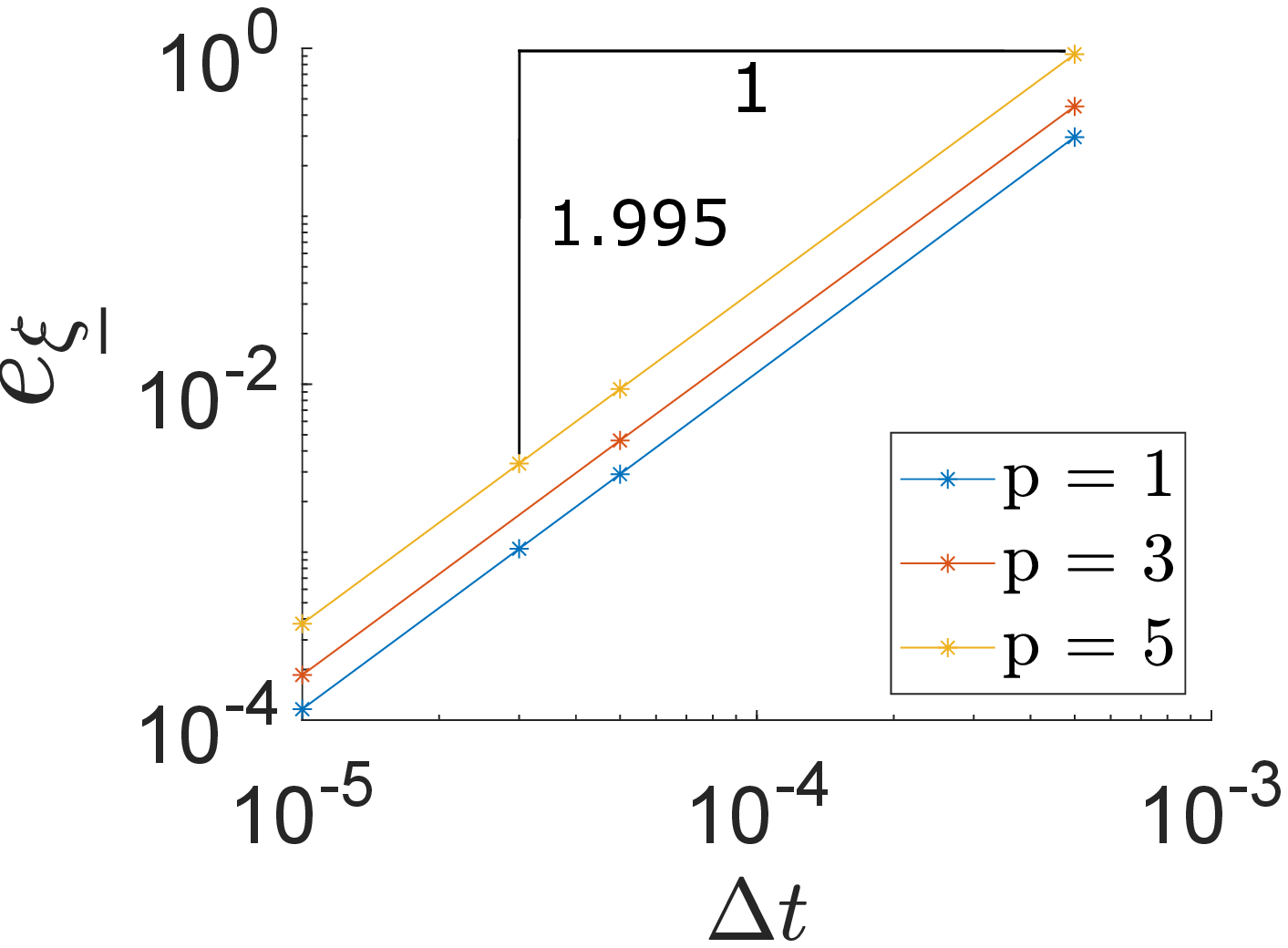}
\end{minipage}
  \hfill
\begin{minipage}{0.315\textwidth}
	\centering
    \includegraphics[width=\textwidth]{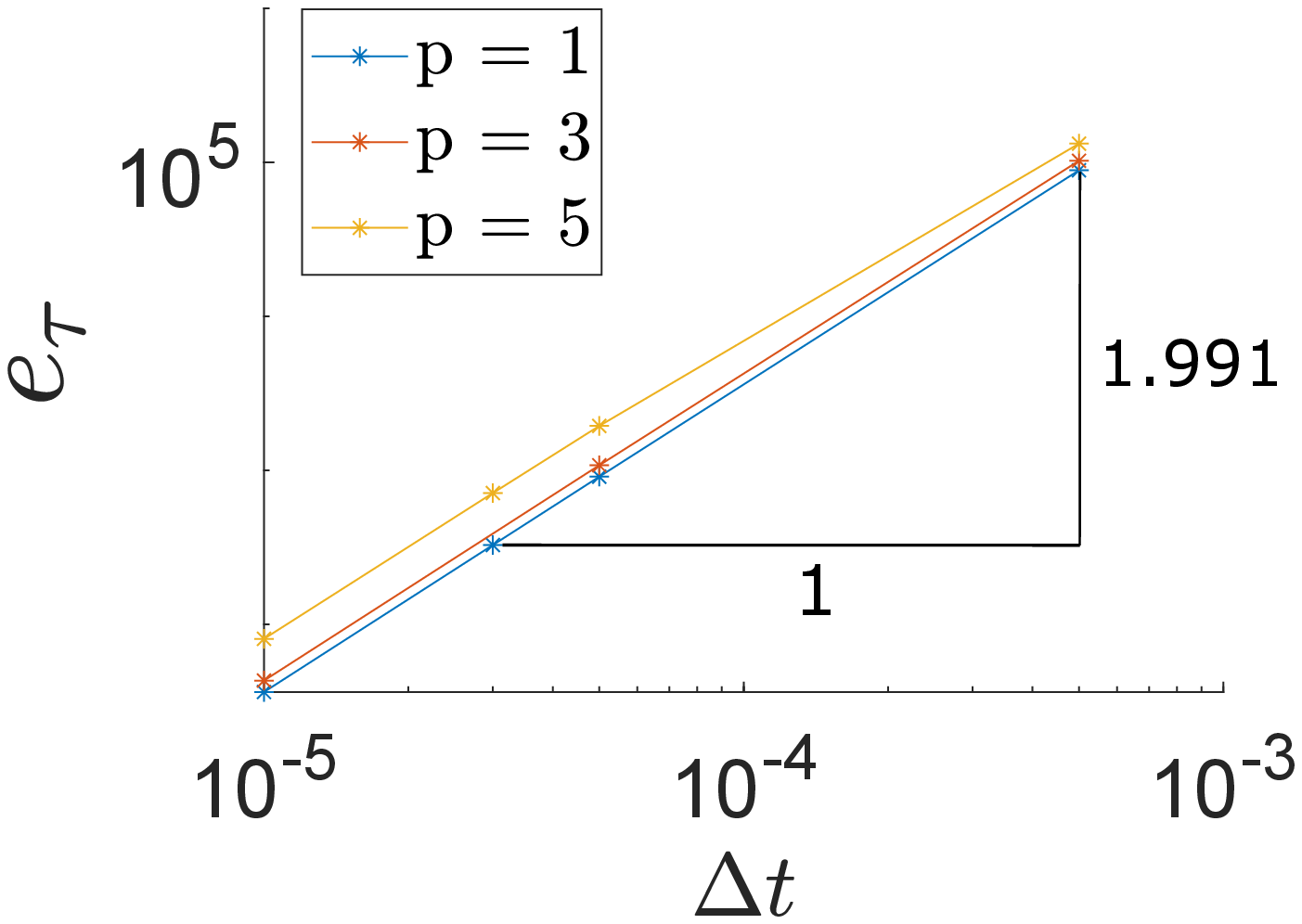}
\end{minipage}\hfill
\begin{minipage}{0.315\textwidth}
	\centering
    \includegraphics[width=\textwidth]{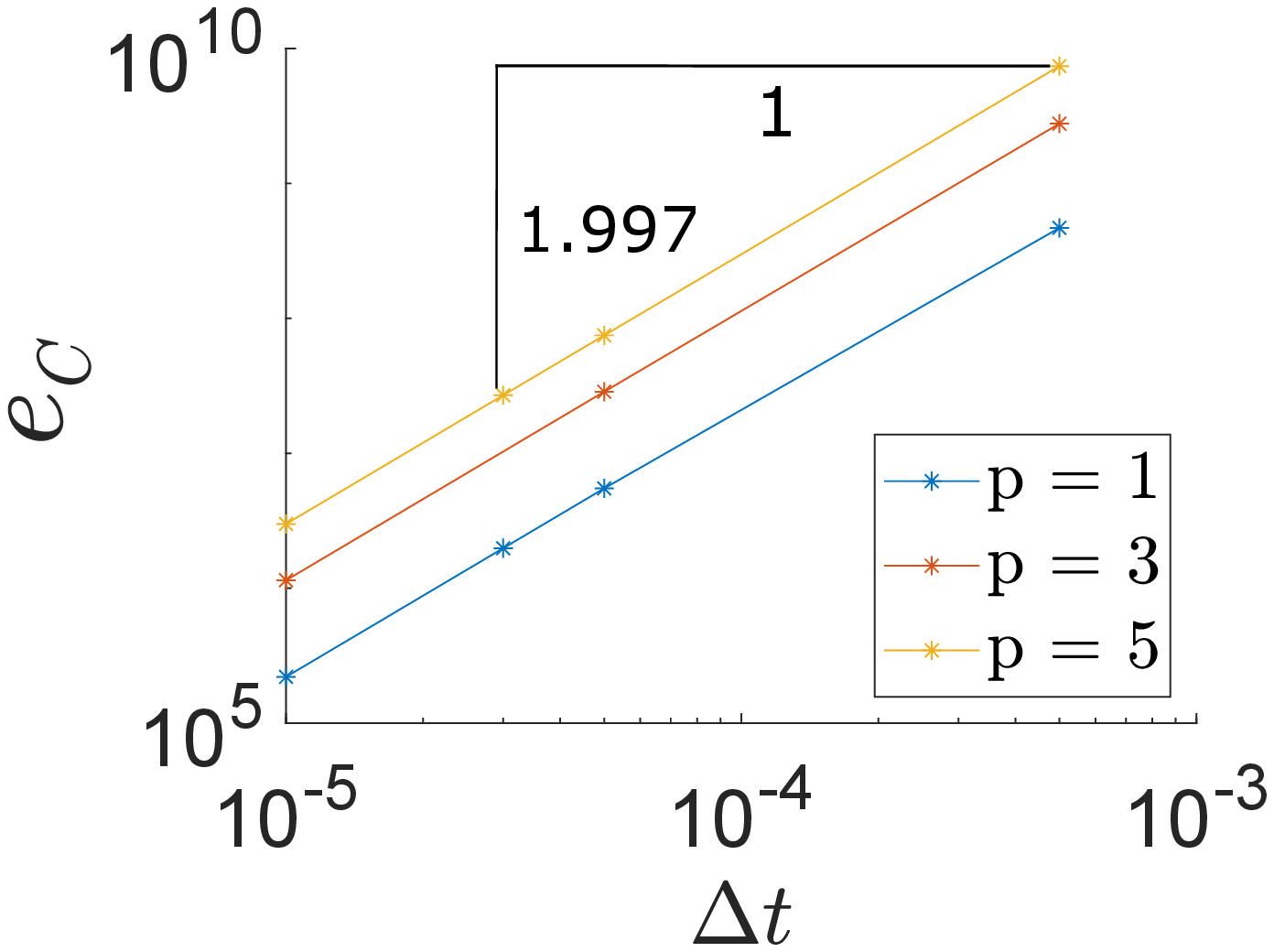}
\end{minipage}
    \caption{Numerical convergence study for the OCP formulation using macro-node error computation for $t_f=0.12s$}
    \label{fig:OrderOCP}
\end{figure}

Next the full multirate DMOC scheme was applied for the solution of the optimal control problem detailed in Equation~(\ref{eq:ProbForm}). The optimization was implemented using the \textit{interior-point-convex} method included in the MATLAB \textit{quadprog} routine and validated against the analytical solution outlined by Turner and Junkins.\cite{Reading494}   
For this purpose a numerical convergence test was performed analogously to the previous section for simulations of length $t_f=0.12s$ and $\theta_{t_f}=20^{\circ}$. The test was performed using an absolute error computation for the optimal cost $C$ and error estimation following the definition in Eq.~(\ref{eq:error}) for the trajectory $\underline{\xi}$ and the control force $\tau$. The results presented in Figure~\ref{fig:OrderOCP} are reasonable as performing the maneuver in such short time interval results in values for the cost function and the control input of order as high as $10^{11}$ and $10^{6}$ respectively. The figures show convergence of order 2 for all considered variables, as expected from the use of midpoint quadrature, thus validating the correct behaviour of the numerical scheme.

\begin{figure}[htb]
\begin{minipage}{0.31\textwidth}
	\centering
	\vspace{10pt}
    \includegraphics[width=\textwidth]{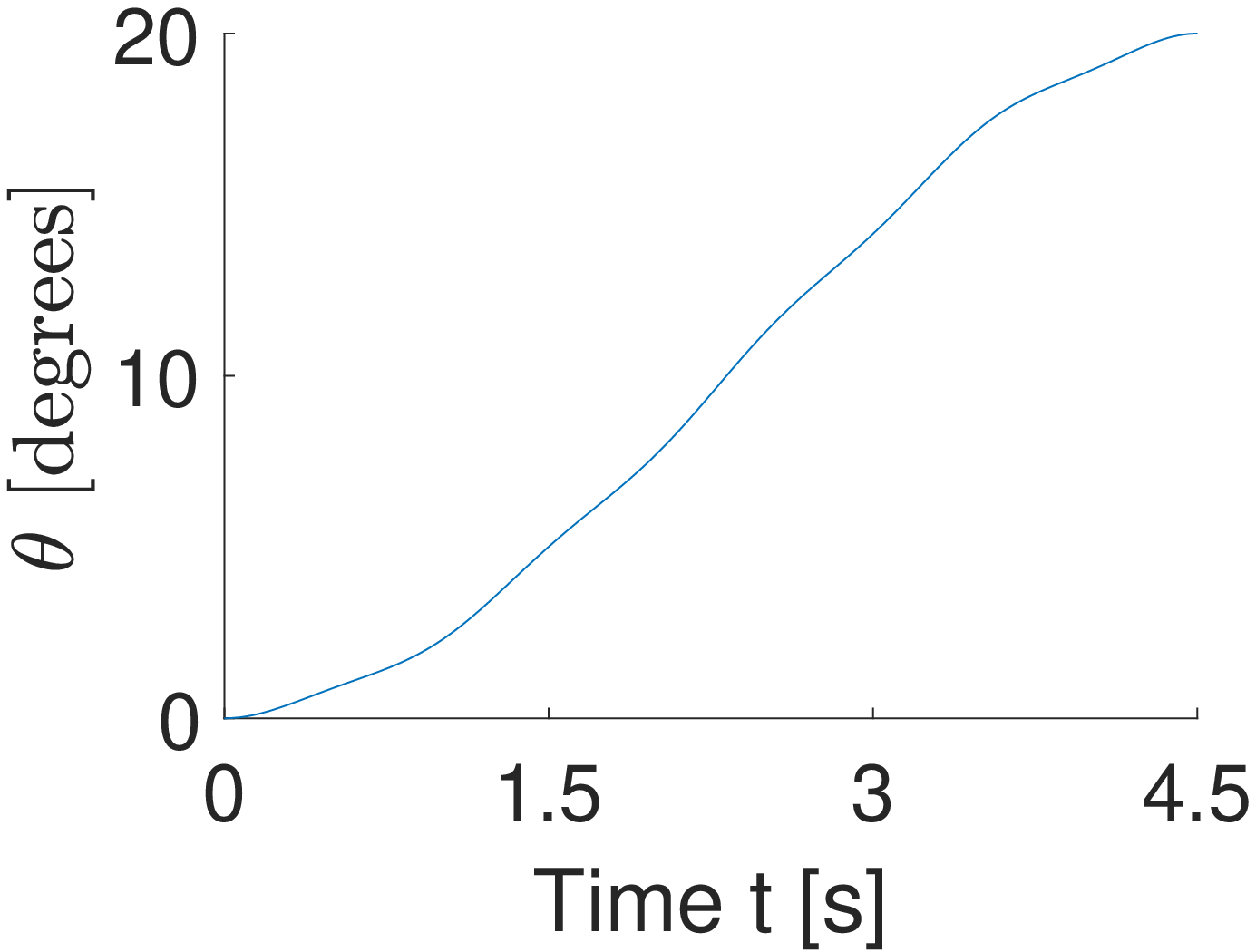}
\end{minipage}\hfill
\begin{minipage}{0.31\textwidth}
	\centering
    \includegraphics[width=\textwidth]{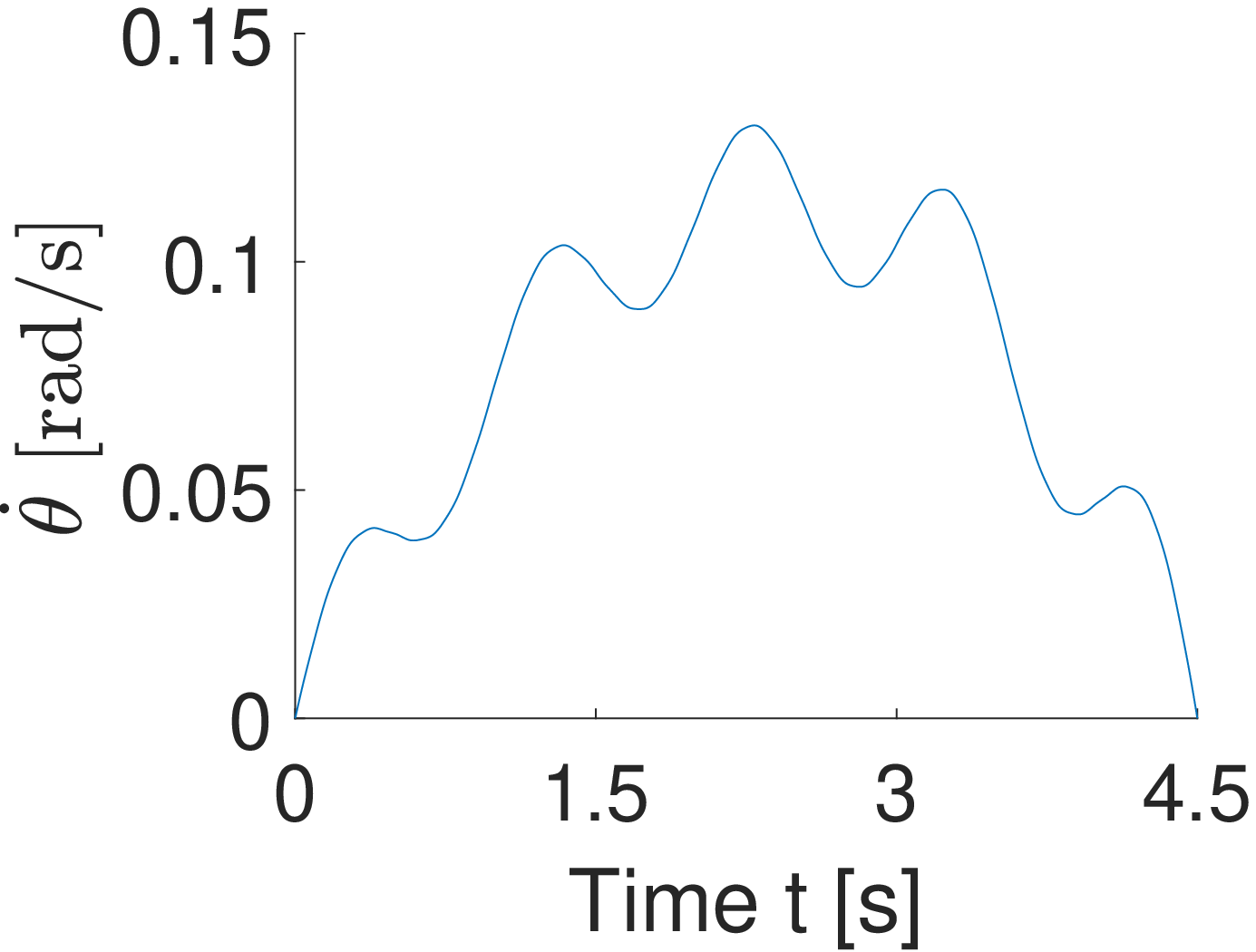}
\end{minipage}
  \hfill
\begin{minipage}{0.31\textwidth}
	\centering
	\vspace{10pt}
    \includegraphics[width=\textwidth]{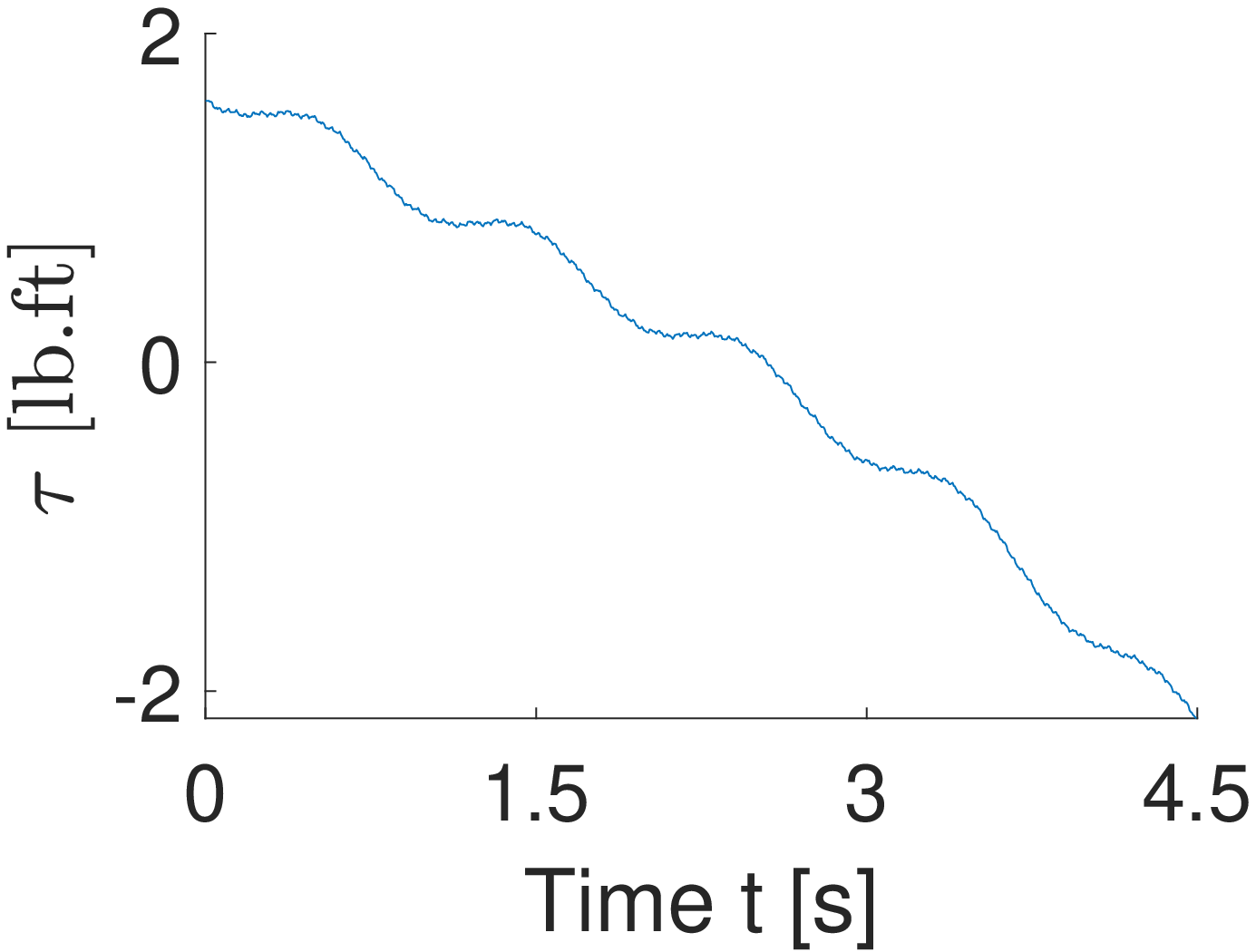}
\end{minipage}

\begin{minipage}{0.31\textwidth}
\vspace{10pt}
	\centering
    \includegraphics[width=\textwidth]{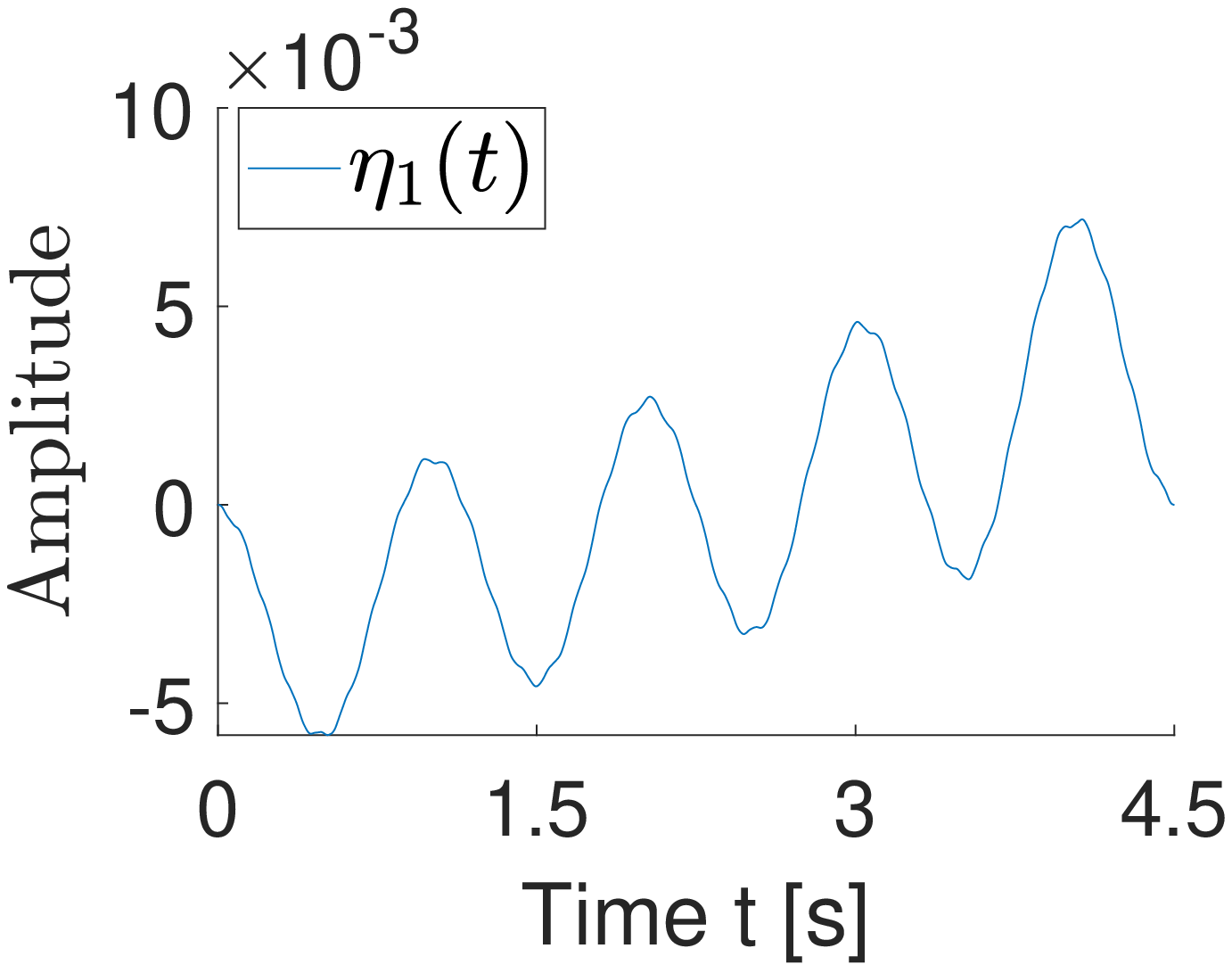}
\end{minipage}\hfill
\begin{minipage}{0.31\textwidth}
	\centering
	\vspace{10pt}
    \includegraphics[width=\textwidth]{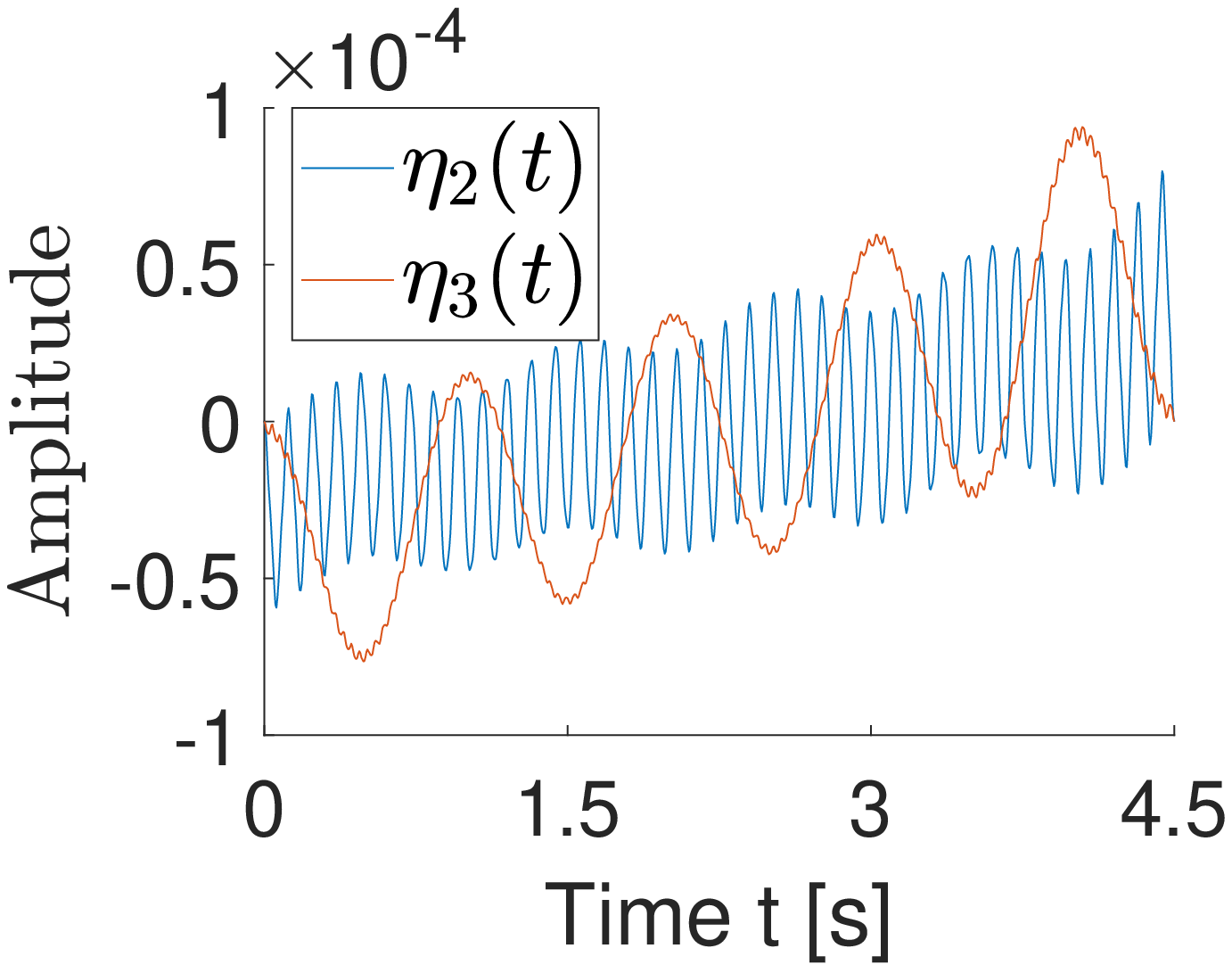}
\end{minipage}
\hfill
\begin{minipage}{0.31\textwidth}
	\centering
	\vspace{10pt}
    \includegraphics[width=\textwidth]{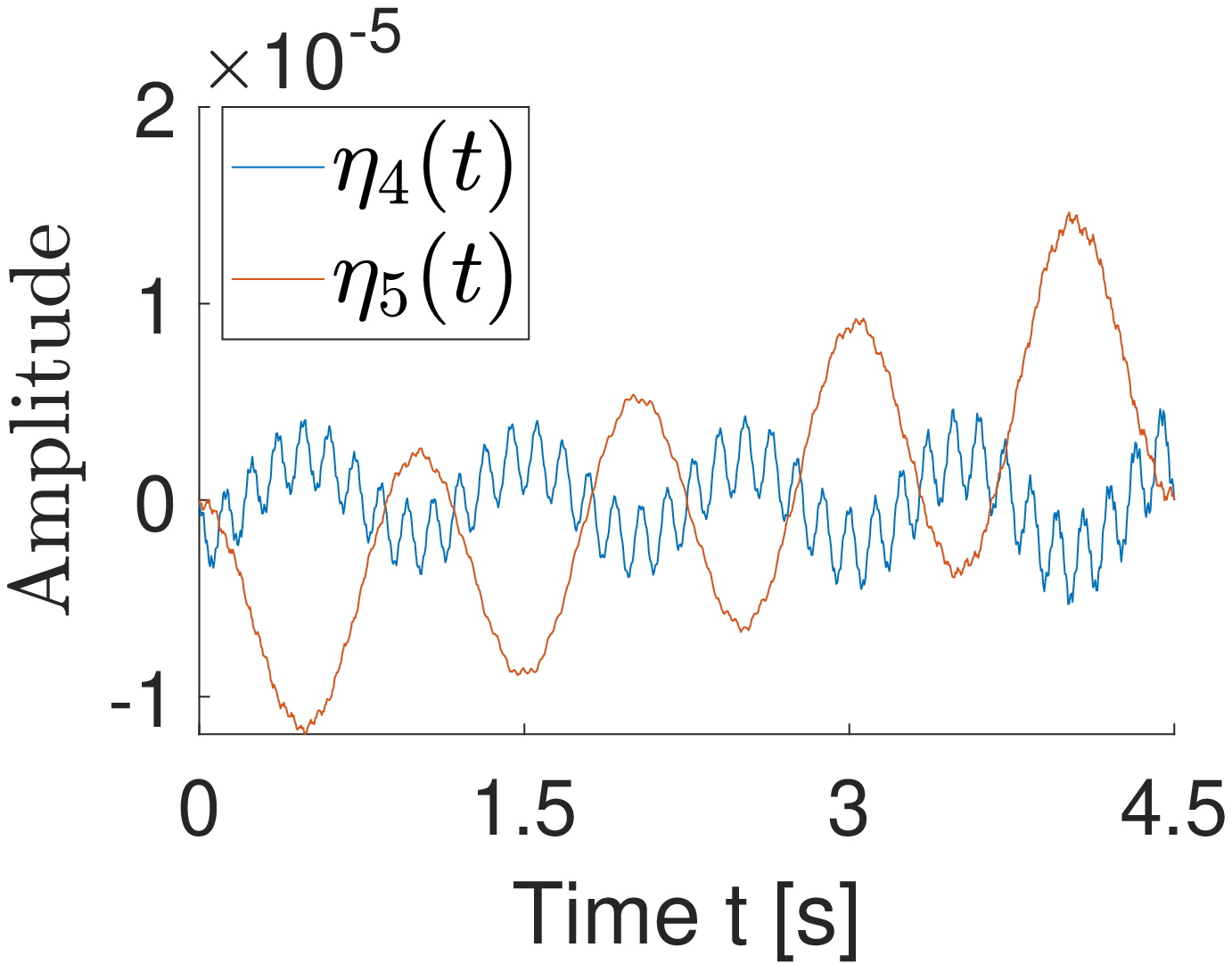}
\end{minipage}
    \caption{OCP solution with Multirate DMOC for $\theta_{t_f}=20^{\circ}$, $t_f=4.5s$, $\Delta t=10^{-3}$ and $p=5$}
    \label{fig:OCPrun}
\end{figure}

\begin{figure}[htb]
\begin{minipage}{0.48\textwidth}
	\centering
    \includegraphics[width=\textwidth]{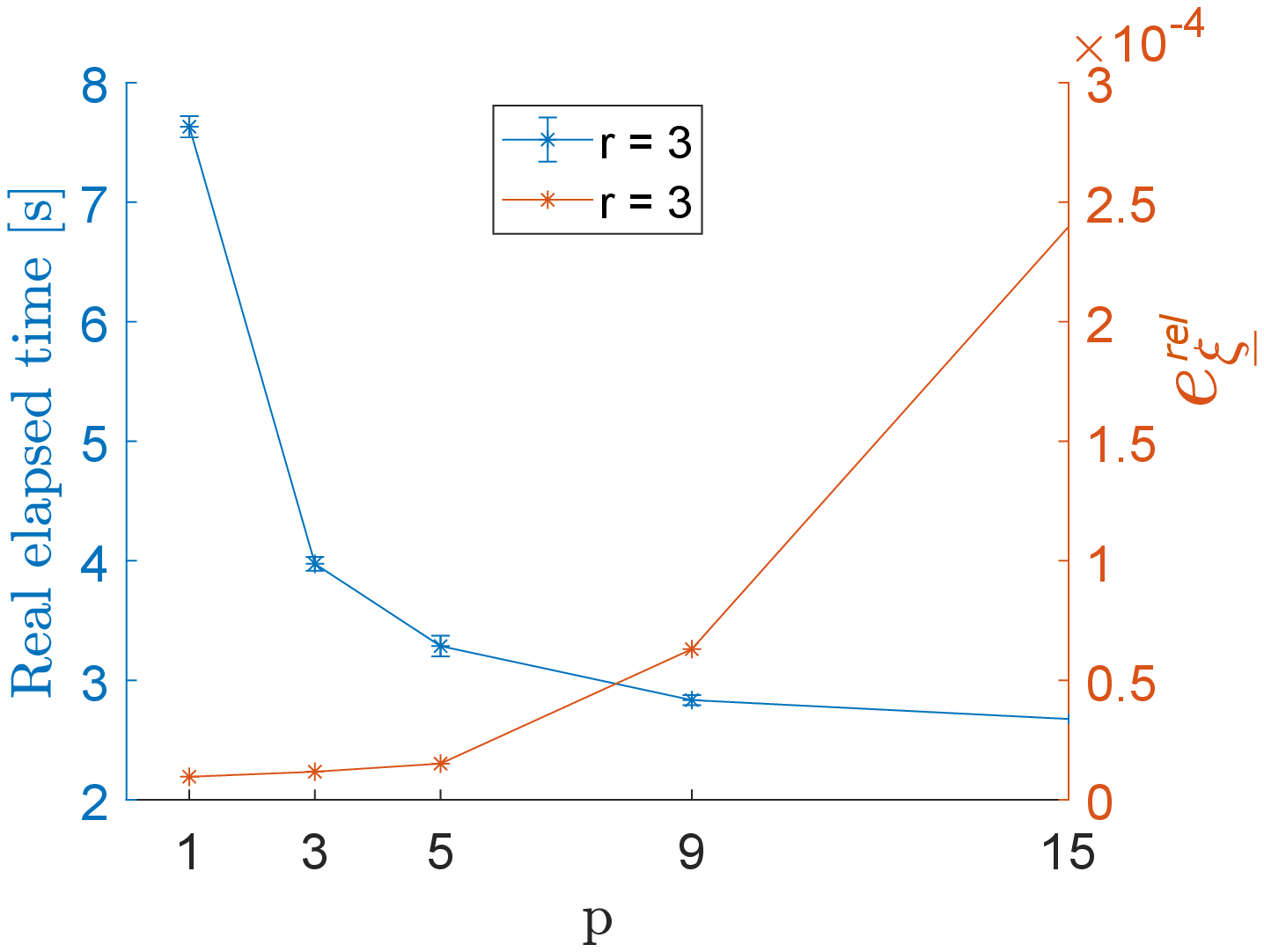}
	\caption{Mean computational time with standard deviation and relative error in $\;\underline{\xi}$ versus $p$ for a constant micro time step of $10^{-3}$, $t_f=4.5s$ and $\theta_{t_f}=20^{\circ}$}
	\label{fig:Tradeoff}
\end{minipage}
  \hfill
\begin{minipage}{0.48\textwidth}
	\centering
	%\vspace{5pt}
    \includegraphics[width=\textwidth]{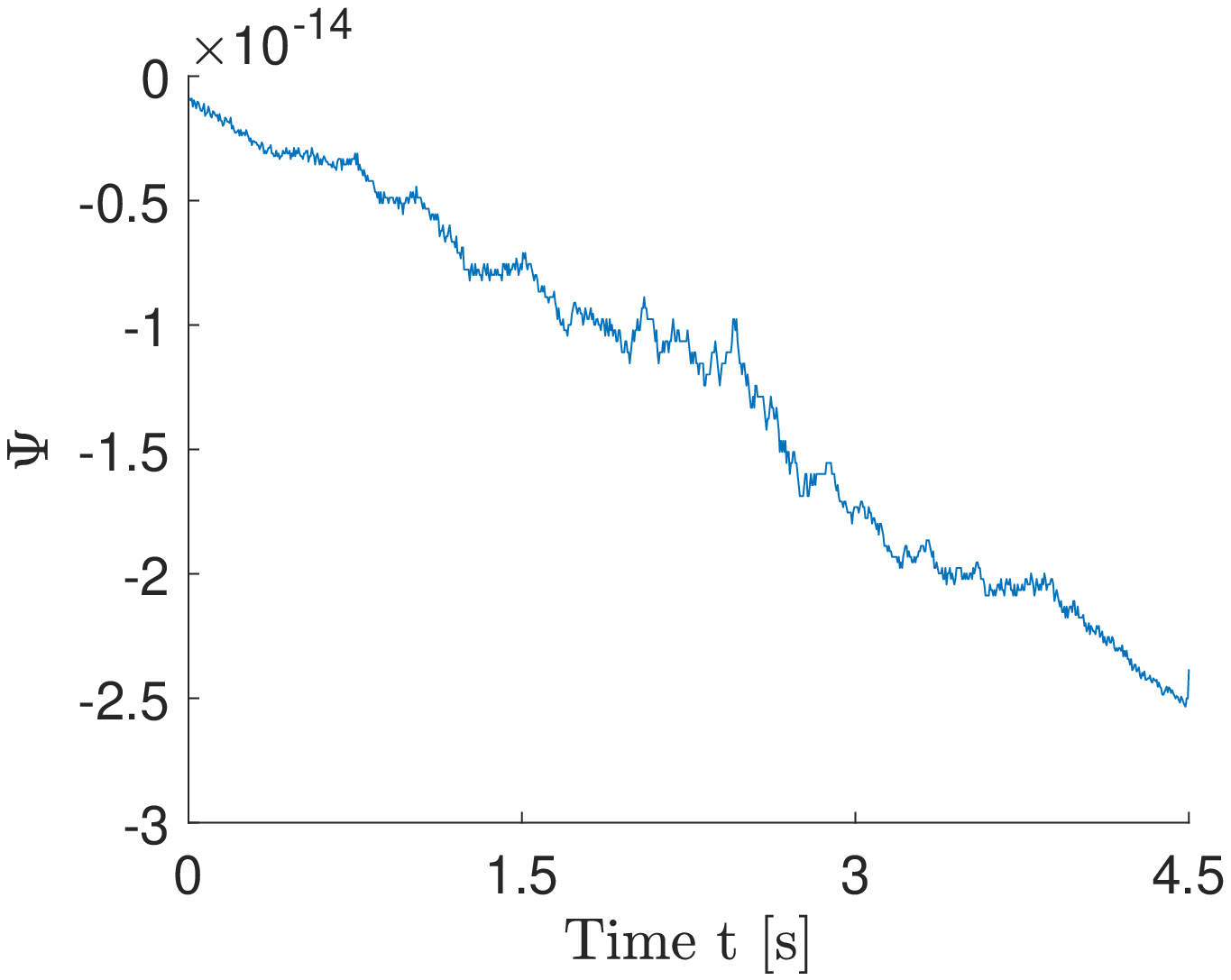}
    \caption{Demonstration of the conservation properties of Multirate DMOC for a simulation with $\Delta t=10^{-3}$, $p=5$, $t_f=4.5s$ and $\theta_{t_f}=20^{\circ}$}
    \label{fig:MomntDisc}
\end{minipage}\hfill
\end{figure}

An example solution of the multirate OCP is presented in Figure \ref{fig:OCPrun} for a rest-to-rest maneuver with $\theta_{t_f}=20^{\circ}$, $t_f=4.5s$, $p=5$ and $\Delta t=10^{-3}$. The results demonstrate the successful execution of the maneuver and vibration suppression in all modes included in the model. %\textcolor{red}{, the achievement of vibration suppression in all modes included in the model} and the keeping/respect of the constraints.
This example solution is also used to demostrate the structure-preserving properties of the Multirate DMOC scheme. In the presence of control forces the difference in total energy and external work is now conserved and the momentum map evolves according to the Discrete Noether's theorem with forcing. \cite{Reading71,Reading57} For the current multirate example the theorem can be shown to be equivalent to the following discrete conservation law
\begin{equation}
\Psi_d^k = p^{\theta}_k-p^{\theta}_0 - \sum_{i=0}^k \sum_{m=0}^{p-1} \Delta t \: \tau_i^{m+1/2} = 0 \; \; \; \; \textrm{for} \; \; k=0,\ldots,n_s-1\;\;\textrm{and}\; \;  \mathbf{E}^T \begin{bmatrix} 
         p_k^{\theta} \\
        \underline{p}_k^{\underline{\eta}}
         \end{bmatrix} = \begin{bmatrix} 
         \underline{p}_k^s \\
         \underline{p}_k^f
         \end{bmatrix}
\end{equation}
based on the linear transformation introduced in Eq.~(\ref{eq:MTransf}). Here $ \underline{p}_k^s = \underline{p}_k^{s+} = \underline{p}_k^{s-}$ and \newline $ \underline{p}_k^f = \underline{p}_k^{f,m+} = \underline{p}_k^{f,m-}$ , which can be derived by replacing Eq.~(\ref{eq:ConjMomentum}) in the discrete equations of motion (Eq.~(\ref{eq:DiscrEqMot1})).
For the example simulation above, Figure~\ref{fig:MomntDisc} demonstrates that the Multirate DMOC method successfully preserves this law and testifies for the structure-preserving properties of the multirate OCP formulation.

\begin{figure}[htb]
  \begin{subfigure}[b]{0.48\textwidth}
    \includegraphics[width=\textwidth]{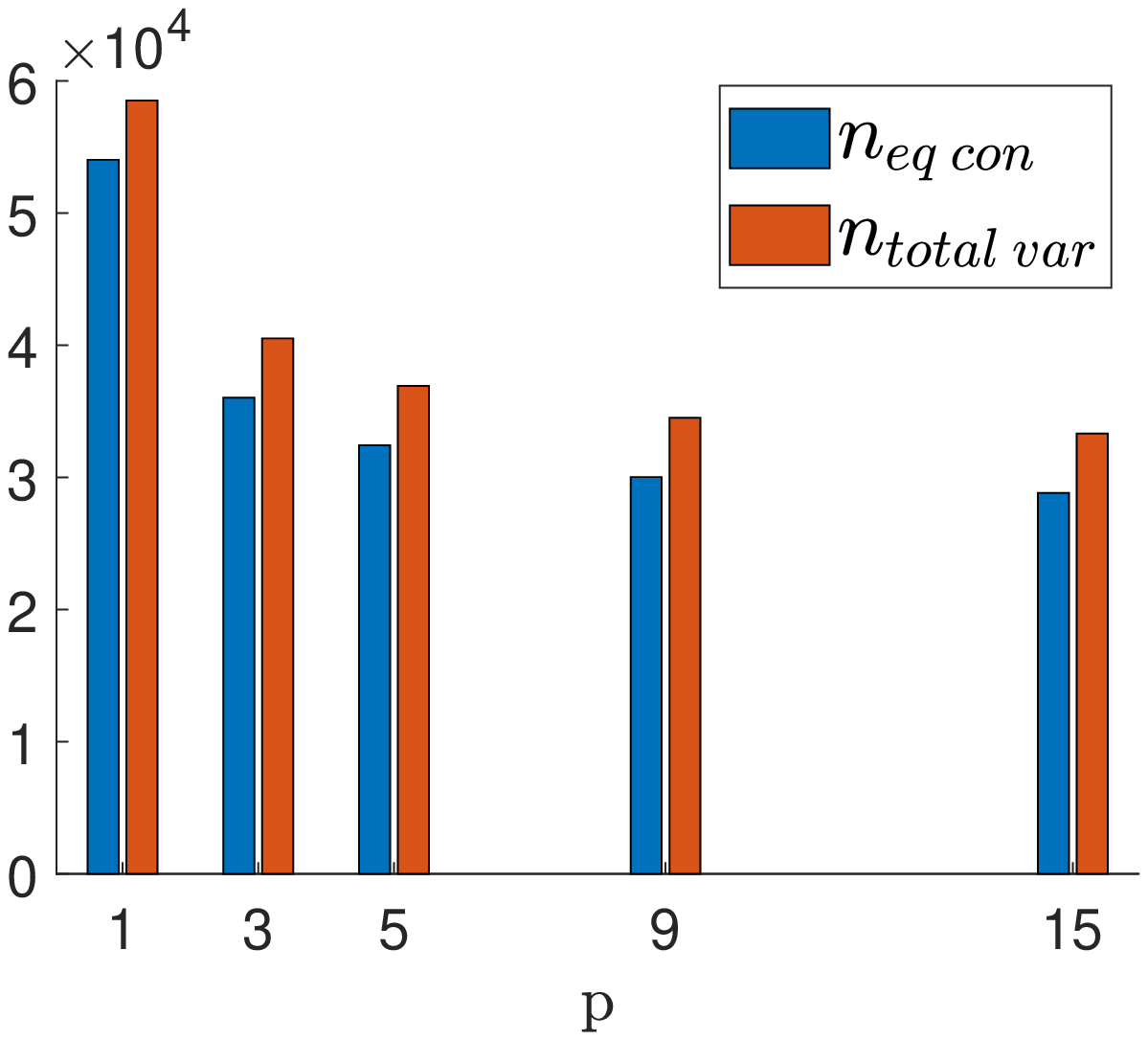}
  \end{subfigure}
  \hfill
  \begin{subfigure}[b]{0.48\textwidth}
    \includegraphics[width=\textwidth]{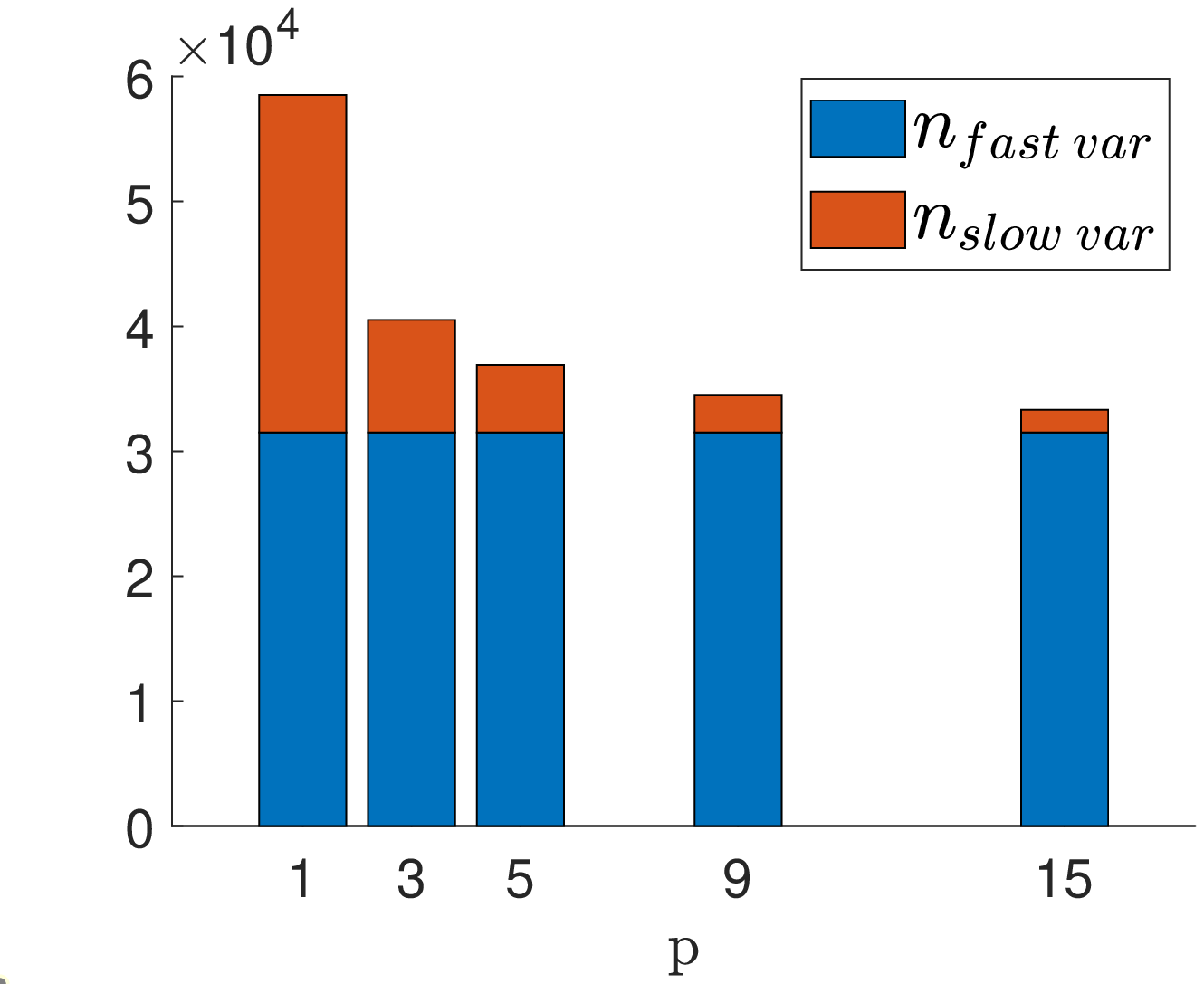}
  \end{subfigure}
    \caption{Size of OCP based on Multirate DMOC for a simulations with $\Delta t=10^{-3}$ and $t_f=4.5s$}
    \label{fig:SizeOCP}
\end{figure}

 Next, the advantages of the multirate OCP scheme were demonstrated in a series of tests comparing the single rate ($p=1$) and the multirate ($p>1$) solutions. Fundamentally, increasing the proportionality $p$ leads to a coarser macro time grid on which the slow variables are resolved lowering the accuracy of their computation. On the other hand, using a coarser macro time grid reduces the number of time nodes on which the slow variables are computed and results in overall reduction in the number of optimization variables $n_{total\, var}$ and equality constraints $n_{eq\,con}$ used to define the optimization problem. This effect is demonstrated in Figure \ref{fig:SizeOCP} for simulations with $\Delta t=10^{-3}$ and $t_f=4.5s$, where $n_{total\:var}\!= n_{slow \,var}+ n_{fast \,var}$, where $n_{slow \,var}$ and $n_{fast \,var}$  denote respectively the number of optimization variables resulting from discretization on the macro and micro grid.

A more detailed look at the trade-off between accuracy of the simulation and computational cost of the problem is presented in Figure \ref{fig:Tradeoff}. In this test a series of simulations is performed keeping the micro time step constant at $10^{-3}$ while varying the macro-micro time step proportionality and recording the resulting variations in both accuracy and real elapsed time for the optimization.   % based on the relative error definition in Eq. (\ref{eq:relError}). % \textcolor{red}{
For this test the relative error is computed using the following definition
 \begin{equation}
	\label{eq:relError}
	e^{rel}_{\underline{\xi}}(\underline{\xi}_{\,d},\underline{\xi})=\frac{ \max_{k\:=\:0,...,n_s} (\| \underline{\xi}_k-\underline{\xi}(t_k)\|_{\infty})}{\max_{k\:=\:0,...,n_s} (\| \underline{\xi}(t_k)\|_{\infty})}
\end{equation}
The computation time %records
in Figure \ref{fig:Tradeoff}  denotes the mean real elapsed time for the execution of \textit{quadprog} obtained from 10 measurements using the MATLAB \textit{tic-toc} routine.  
From the figure it can be seen that  by adopting the proposed method it is possible to reduce computational costs significantly for small sacrifices in accuracy. For example in our case the elapsed runtime is more than halved between $p = 1$ to $p = 5$ for a negligible increase in the relative error $e^{rel}_{\underline{\xi}}$ from $9.632\times10^{-6}$ to $1.514\times10^{-5}$. Furthermore, unlike the real elapsed time curve for the numerical integrator, the real elapsed time for the \textit{quadprog} function monotonically decreases. However, at greater values of $p$ we  observe diminishing returns as the error increases rapidly for smaller computational time reductions. Thus the multirate scheme allows for great flexibility, allowing the designer to decide upon the trade-off between computational cost and simulation fidelity. 

The trade-off can further be optimized by investigation of the optimal combination of proportionality $p$ and number of variables $r$ discretized on the macro grid as slow variables. For this implementation of multirate DMOC the dependence of the size of the number of optimization variables can be expressed as follows
 \begin{equation}
	\label{eq:ten}
	n_{total\:var} = n_{slow\: var}+n_{fast\: var}
\end{equation}
 \begin{equation}
	\label{eq:eleven}
	n_{slow\: var}(p,r,N,t_f,\Delta t) = 2\:r\,\Big(\frac{t_f}{p\:\Delta t}+1\Big)
\end{equation}
 \begin{equation}
	\label{eq:fourteen}
	n_{fast\: var}(r,N,t_f,\Delta t) = 2\:(N+1-r)\Big(\frac{t_f}{\Delta t}+1\Big)+\frac{t_f}{\Delta t}
\end{equation}

\begin{figure}[htb]
\begin{minipage}{0.48\textwidth}
\vspace{10pt}
	\centering
    \includegraphics[width=\textwidth]{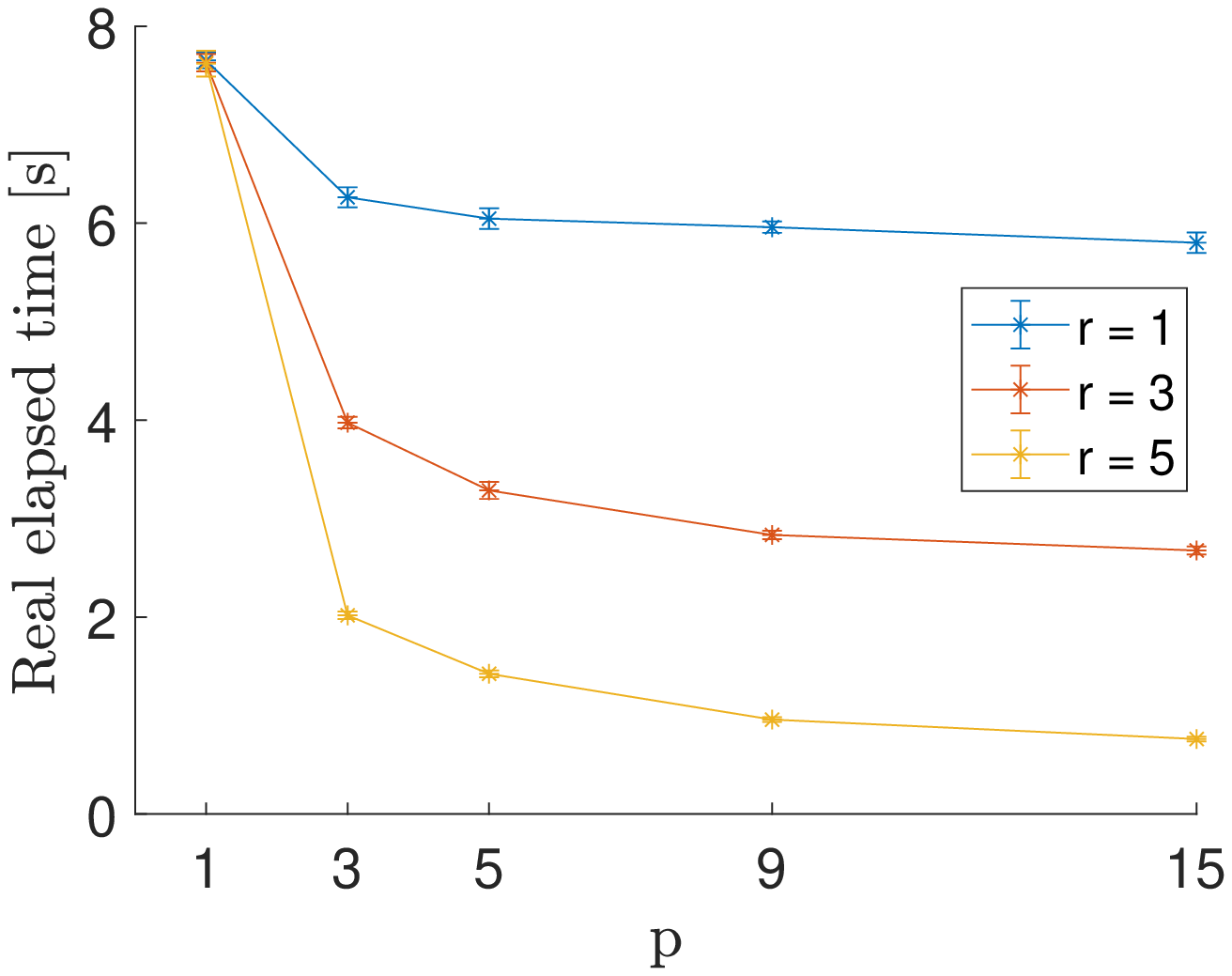}
    \caption{Mean computational time with standard deviation versus $p$ for a constant
micro time step of $10^{-3}$, $t_f=4.5s$ and $\theta_{t_f}=20^{\circ}$}
    \label{fig:AccCPUtime}
\end{minipage}
  \hfill
\begin{minipage}{0.48\textwidth}
	\centering
	\vspace{10pt}
    \includegraphics[width=\textwidth]{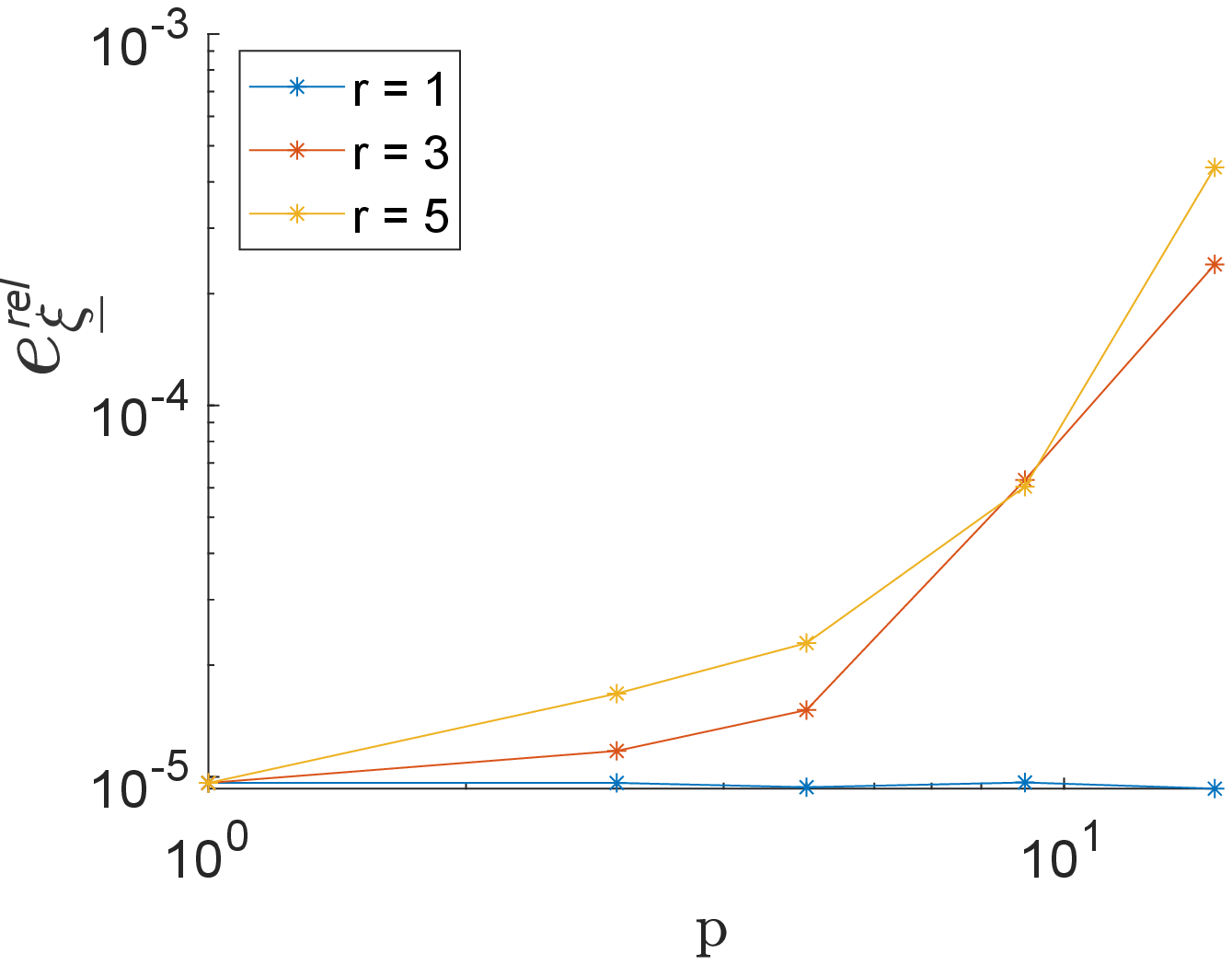}
    \caption{Relative error in $\underline{\xi}$ versus $p$ for a constant
micro time step of $10^{-3}$, $t_f=4.5s$ and $\theta_{t_f}=20^{\circ}$}
    \label{fig:Acc}
\end{minipage}\hfill
\end{figure}

As presented in Figures~\ref{fig:AccCPUtime} and ~\ref{fig:Acc} increasing $r$ and thus reducing the variables discretized on the micro grid allows for even larger reductions in computational cost. The choices in $p$ and $r$ provide freedom to the practitioner to tailor the method to the time-scales present in the problem and allow one to obtain a high fidelity solutions at reduced computational cost.

\newpage
\vspace*{-35pt}
\section{Conclusion}
In this work we propose the use of a multirate variational optimal control scheme known as Multirate DMOC for simultaneous attitude and vibration control of flexible spacecraft. These structures experience dynamics on multiple time scales and thus present competing challenges for the numerical integrator and the optimal control solver. The use of small time steps ensures the correct resolution of fast evolving dynamics, but results in unnecessary computational cost for the approximation of the slow subsystem. For this purpose the proposed method separates the system into a slow and a fast subsystem and respectively discretizes them on a macro and micro time grid. The multirate equations of motion and the necessary optimality conditions for the OCP are obtained by direct discretization of the variational principle. The resulting structure-preserving time-stepping equations serve as equality constraints for the optimization problem and allow for a discrete OCP formulation, which inherits the conservation properties of the continuous-time model.

To demonstrate the advantages of multirate DMOC we construct a general linear model of a flexible spacecraft and formulate an optimal control problem to perform a single-axis rotation while leaving all considered modes of vibration quiescent at the end of the maneuver. Comprehensive investigations for this example system demonstrate the numerical convergence and conservation properties of both the numerical integrator and the full multirate OCP scheme. Ultimately, the multirate discretization leads to a reduction in the number of optimization variables and equality constraints whilst providing a sparse structure for the constraint Jacobian. In a series of simulations it is shown that by tailoring the macro-micro time step proportionality $p$ one can achieve significant reductions in computation cost for a negligible penalty in accuracy. For the specific example the real elapsed time for the optimization is more than halved for $p$=$5$ for a negligible penalty in accuracy.
Furthermore it is demonstrated that in systems with several slow motions, the number of generalized coordinates discretized on the macro scale can straightforwardly be extended within the proposed method allowing for greater reductions in the computational cost. Thus the multirate formulation provides freedom to the designer to customize the method to the time-scales present in the problem, allowing for a high fidelity solution at a reduced computational cost. Future work will investigate procedures for finding optimal $p$ and $r$ values and examine %study 
the performance of Multirate DMOC for flexible spacecraft models including kinematic nonlinearities and dissipation effects.

\bibliographystyle{unsrt}   % Number the references.
\bibliography{references}   % Use references.bib to resolve the labels.

\end{document}